\documentclass{article}
\usepackage{graphicx} %
\usepackage{mycommands}
\usepackage{xcolor}
\usepackage{natbib}
\usepackage{hyperref}
\usepackage{slashbox}
\usepackage{enumitem}
\usepackage{tikz}
\usepackage{caption}
\usepackage{subcaption}
\usepackage{booktabs} 
\usepackage{slashbox}

\usepackage{amsmath}
\usepackage{tabularx}
\usepackage{mathtools}

\usepackage{lmodern}
\usepackage{blkarray}
\usepackage{bigstrut}
\usepackage{mathtools}
\usepackage{orcidlink}
\usepackage{authblk}
\usepackage{svg}

\title{Quantifying and testing dependence to categorical variables}
\author[1]{Siegfried H\"ormann\thanks{\href{mailto:shoermann@tugraz.at}{shoermann@tugraz.at, \orcidlink{0000-0001-5878-2840}}}}
\author[1]{Daniel Strenger-Galvis\thanks{\href{mailto:daniel.strenger@tugraz.at}{daniel.strenger@tugraz.at, \orcidlink{0009-0002-2465-604X}}}}
\affil[1]{Institute of Statistics, Graz University of Technology}

\date{\today}

\setcitestyle{round}

\begin{document}

\maketitle
\begin{abstract}
    We suggest a dependence coefficient between a categorical variable and some general variable taking values in a metric space. We derive important theoretical properties and study the large sample behaviour of our suggested estimator. Moreover, we develop an independence test which has an asymptotic $\chi^2$-distribution if the variables are independent and prove that this test is consistent against any violation of independence. The test is also applicable to the classical~$K$-sample problem with possibly high- or infinite-dimensional distributions. We discuss some extensions, including a variant of the coefficient for measuring conditional dependence.
\end{abstract}


\section{Introduction}\label{sec:introduction}
Let~$Y$ be a categorical random variable with~$K$ levels, and let~$X$ be a covariate taking values in a separable metric space. The aim of this paper is to introduce a dependence coefficient that quantifies the strength of association between~$Y$ and~$X$, and subsequently to develop a statistical test for the null hypothesis
$$
\mathcal{H}_0\colon \text{$Y$ and~$X$ are independent.}
$$

The problem of measuring dependence between random variables has experienced a significant upswing in the past few years, particularly following the influential work of \cite{chatterjee_new_2021}, who addressed this question in the setting where both~$Y$ and~$X$ are real-valued. He proposed the dependence coefficient
\begin{equation}\label{eq:xi-definition}
    \xi(X,Y) = \frac{\int \Var(P(Y\geq t|X))dF(t)}{\int \Var(1\lbrace Y\geq t\rbrace)dF(t)},
\end{equation}
with~$F(t)=P(Y\leq t)$. Unlike Pearson’s correlation coefficient, Spearman’s~$\rho$ or Kendall’s~$\tau$, which are often used in applications, the coefficient~$\xi$ is characterizing independence. In particular, the following key properties can be shown to hold:
\begin{enumerate}[label = (M\arabic*)]
\item\label{it:m1}~$0\leq \xi(X,Y)\leq 1$,
    \item\label{it:m2}~$\xi(X,Y)=0$ if and only if~$X$ and~$Y$ are independent,
    \item\label{it:m3}~$\xi(X,Y)=1$ if and only if there is a measurable function~$f:\R\to\R$ such that~${Y=f(X)}$ almost surely.
\end{enumerate}

Recent related contributions include \cite{kroll_asymptotic_2025}, who shows that~$\xi$ is asymptotically normal whenever~$Y$ is not constant, \cite{dette_simple_2025}, who provide an estimator for the asymptotic variance (under dependence or independence) and \cite{azadkia_new_2025}, who introduce a related coefficient that appears to be more powerful for detecting certain types dependencies.
An overview of several other dependence coefficients that have been previously studied in the literature can also be found in \cite{chatterjee_new_2021} and \cite{chatterjee_survey_2024}. For the sake of brevity, we do not provide yet another survey.

In its population setting the coefficient~$\xi$ can be immediately extended to a setup where~$X$ takes values in more general spaces. \cite{azadkia_simple_2021} have considered multivariate covariates and---motivated by functional data---\cite{hormann_azadkiachatterjees_2024} have studied infinite-dimensional~$X$. While there is little conceptual difference in the population version of these extensions, substantial technical challenges arise in the estimation and the derivation of large-sample properties of the coefficient.

\par The previously cited papers crucially use the fact that $Y$ is real valued and exploit the ordering of $\mathbb{R}$. Alternative methods that are capable of handling multivariate~$Y$ are e.g.\  distance correlation \citep{szekely_measuring_2007} or the HHG-test \citep{heller_consistent_2012}. More recently, \cite{deb_measuring_2020} proposed a new class of non-parametric dependence measures that include~$\xi$ as a special case and even allow both variables $X$ and $Y$ to take values in general topological spaces.

\par In this paper, we turn to the presumably more specialized situation, where~$Y$ is a categorical random variable.  \cite{chatterjee_new_2021} noted that the dependence measure~$\xi$ can also be applied in this case, \emph{“\ldots by converting the categorical variables to integer-valued variables in any arbitrary way.”} Properties~{\ref{it:m1}--\ref{it:m3}} continue to hold.  However, although~$\xi$ is invariant under monotonic transformations, this is generally not true for permutations of the values of~$Y$. Thus, if the categories of~$Y$ are encoded as the integers~${1, \ldots, K}$, the specific ordering of this assignment will generally affect the value of the measure. The following result shows that such effects can be arbitrarily strong: 
\begin{Proposition}\label{pr:1} 
For some categorical variable~$Y$ let~$\mathcal{L}_Y$ be the class of bijections between the~$K$ levels of~$Y$ and the integers~$\{1,\ldots, K\}$.
For any~$\varepsilon>0$ there exists a categorical variable~$Y$ and a random variable~$X$ such that
$$
\max_{\ell\in\mathcal{L}_Y}\xi(X,\ell(Y))-\min_{\ell\in\mathcal{L}_Y}\xi(X,\ell(Y))>1-\varepsilon.
$$
\end{Proposition}
Proposition~\ref{pr:1} highlights that the interpretation of a dependence measure for a categorical variable
$Y$ is misleading if it is not invariant under the transformations~$\ell\in\mathcal{L}_Y$. Similar effects arise if we adapt other prominent dependence measures, like for example \textit{distance correlation} \citep{szekely_measuring_2007}.

\par The relationship between $X$ and $Y$ can be recast in terms of conditional distributions of $X$. Independence corresponds to the classical $K$-sample problem: if $X$ and $Y$ are independent, then the conditional distributions of $X$ given different values of $Y$ coincide; dependence arises when at least one category induces a different distribution. At the opposite extreme, functional dependence is equivalent to mutual disjointness of the conditional distributions: each value of $Y$ corresponds to a subset of $X$’s image space with no overlap across categories.

\par This perspective connects our study to established testing frameworks. In low-dimensional settings, the $K$-sample problem has been extensively addressed by rank- and distance-based procedures such as the Kruskal–Wallis \citep{kruskal_use_1952}, Kolmogorov–Smirnov \citep{kolmogorov_sulla_1933}, and Cramér–von Mises \citep{anderson_distribution_1962} tests, as well as by divergence-based approaches from information theory (e.g., Kullback–Leibler or Jensen–Shannon divergence). These methods, however, do not extend well to high- or infinite-dimensional $X$. For such cases, alternative approaches have been proposed, e.g., by \cite{schilling_multivariate_1986} and \cite{hall_permutation_2002}. Note that both are limited to the two-sample setting (${K=2}$). For the more general~$K>2$ case in high dimensions, \cite{balogoun_k-sample_2022} introduced a variant of the Maximum Mean Discrepancy. Among these high-dimensional methods, Schilling’s test is the only one with a tractable asymptotic null distribution; all others require computationally intensive resampling. Importantly, while these procedures address independence testing, none provide a quantitative measure of dependence (or, equivalently, separation of conditional distributions) that satisfies properties~\ref{it:m1}–\ref{it:m3}.

\par Measures that allow for quantification of dependence and are invariant under~$\mathcal L_Y$ have been recently proposed by \cite{dang_new_2021} (\textit{Gini distance correlation}) and  \cite{liu_measuring_2025} (\textit{projection correlation}). Those papers, however, also discuss a $K$-sample type context where a scalar or multivariate variable~$X$ is viewed as the response and where the categorical variable~$Y$ is the predictor. While these measures can be used to characterize independence (both satisfy properties~\ref{it:m1} and~\ref{it:m2}), condition~\ref{it:m3} is no longer valid. Instead, the respective coefficients are equal to 1 if and only if the conditional distribution of~$X$ given the category~$Y$ is a single point mass distribution. Note that this implicates that~$X$ needs to be a categorical variable as well.

\par The objective of this paper is to obtain a permutation invariant dependence measure when~$Y$ is a categorical response and~$X$ is the covariate. Moreover, we allow~$X$ to take values in a general separable metric space which may be infinite dimensional. We subsequently derive a statistical test for independence. Although most of the paper is guided by the view of the categorical variable being the response, this test is also a contribution to the important~$K$-sample problem. An advantage compared to several approaches in this greater context is that we obtain a pivotal limiting distribution for our test statistic. Hence, we can easily obtain critical values and do not require computationally intensive permutation tests.

\par The remainder of the paper is organized as follows. In Section~\ref{sec:measure}, we introduce our proposed dependence measure in both, population and sample forms, along with their theoretical properties. Section~\ref{sec:testing} establishes the asymptotic distribution of the test statistic under the null hypothesis of independence. In Section~\ref{sec:extensions} we discuss some directions for extending our method. In particular, we provide a coefficient to measure conditional dependence. Section~\ref{sec:illustration} presents numerical experiments that illustrate the empirical performance of our methods. Finally, Section~\ref{sec:proofs} contains the proofs of our theoretical results.

\par We conclude this section by summarising some general notation that we are going to use in the sequel. Paper-specific notation is introduced at first appearance.  Let \( X \) be a random variable and \( P \) a probability measure. The push-forward measure is denoted by~\( P_X(A) := P(X \in A) \). The chi-squared distribution with \( k \) degrees of freedom is written~\( \chi^2(k) \), and~\( \Phi(z) \) denotes the CDF of the standard normal distribution. For a random variable~\( Z \), the \( \beta \)-quantile is denoted \( q_\beta(Z) \). We use \( \convd \) and \( \convP \) for convergence in distribution and convergence in probability, respectively, while $\stackrel{\mathcal{D}}{=}$ is used to indicate that two variables have identical distribution. The boundary of a set \( A \) is written~\( \partial A \). The Kronecker product is denoted \( \otimes \), and \( \mathrm{vec} \) refers to the vectorization operator that stacks the columns of a matrix into a single vector. We let $\|\cdot\|_F$ denote the Frobenius norm of a matrix. Finally,~$\mathrm{deg}(g)$ refers to the maximum degree of a  vertex within a graph $g$.

\section{The dependence measure}\label{sec:measure}
Assume that~$Y$ has levels~$\{1,\ldots, K\}$ and that~$X$ is a variable taking values in a separable metric space~$(\mathcal{H},d)$. The target is to construct a dependence measure~$\psi(X,Y)$, such that \ref{it:m1}--\ref{it:m3} hold and~$\psi(X,Y)=\psi(X,\ell(Y))$ for all~$\ell \in\mathcal{L}_Y$. 
 
 \subsection{Coupling approach}
 One approach for motivating the coefficient $\xi$ is based on coupling. A coupled version $Y'$ of $Y$ is defined in such a way, that the dependence between $Y$ and $X$ can be related to the dependence between $Y$ and $Y'$. In particular, if $Y$ has a continuous distribution then~$(F(Y),F(Y'))$ has uniform marginals and dependence can be related to the corresponding copula. See e.g.\ \cite{fuchs_quantifying_2024}. We also take this route and reduce our problem to measuring dependence between two categorical variables. To this end, we first need the following representation of~$Y$, which is easy to show.
 \begin{Lemma}\label{lem:coupling}
 Let~$U$ be uniformly distributed on the interval~$(0,1)$ and independent of~$X$. Set~${Q_k(X)=P(Y=k|X)}$ and define the variable 
 \begin{equation}\label{eq:y-representation}
 \widetilde Y = f(X,U)\vcentcolon=\min\left\lbrace k\in \{1,\ldots, K\} \colon\sum_{j=1}^k Q_j(X)>U\right\rbrace.
 \end{equation}
Then~$(\widetilde Y,X)$ and~$(Y,X)$ are identically distributed. 
 \end{Lemma}
By Lemma~\ref{lem:coupling} we may assume that~${Y=f(X,U)}$.
Now define~${Y'=f(X,U')}$, where~$U'\stackrel{\mathcal{D}}{=} U$ is independent of~$X$ and~$U$, yielding an identically distributed couple~$(Y,Y')$. The following observation is crucial for our approach:

\begin{Lemma}\label{lem:x-ind-y}  
        (i)~The variables~$X$ and~$Y$ are independent if and only if~$Y$ and~$Y'$ are independent.
        (ii)~There exists~$f\colon \mathcal{H}\to \{1,\ldots, K\}$ such that~$Y=f(X)$, if and only if there is a function~$\ell\in\mathcal{L}_Y$ such that~$Y=\ell(Y')$.
\end{Lemma}

Lemma~\ref{lem:x-ind-y} implies that the dependence between~$X$ and~$Y$ translates into the dependence between the two categorical random variables~$Y$ and~$Y'$. For those variables in turn, we can use existing dependence coefficients. In this paper we employ \emph{Cram\'er's V} (see \cite{cramer_mathematical_1999}). Let
\begin{equation*}    p_{i,j}=P(Y=i,Y'=j),\qquad p_{i}=P(Y=i)\qquad\text{and}\qquad q_{ j}=P(Y'=j).
\end{equation*}
Cram\'er's~$V$ is then given as
\begin{equation*}
    V(Y,Y')=\frac{1}{K-1}\sum_{i,j}\frac{(p_{i,j}-p_{i}q_{j})^2}{p_{i}q_{j}}.
\end{equation*}

In view of Lemma~\ref{lem:x-ind-y} it makes sense to define a dependence coefficient~$\psi$ between a general variable~$X$ and a categorical variable~$Y$ as 
$$
\psi(X,Y)\vcentcolon=V(Y,Y').
$$
\begin{Remark}
    Note that~$p_{j}=q_{j}$. Hence, in the population version we will replace the $q_j$'s by~$p_j$'s.
\end{Remark}

\begin{Remark}
     The statement of our Lemma~\ref{lem:x-ind-y} is analogue to Theorem~1 of \cite{fuchs_quantifying_2024}, where it is formulated for continuous variables in terms of copulas.
\end{Remark}

\begin{Remark}\label{rem:frobenius}
    It can be readily verified that
   ~$$
    p_{i,j}-p_{i}p_{j}=\mathrm{Cov}(Q_i(X),Q_j(X)).
   ~$$
    Thus, setting~$Q(X)=(Q_1(X),\ldots, Q_K(X))'$ and defining~$D=\mathrm{diag}(p_{1},\ldots, p_{K})$, we may alternatively write
    \begin{align}
    (K-1)\psi(X,Y)&=\|D^{-1/2}\,\mathrm{Var}(Q(X))\,D^{-1/2}\|_F^2\label{eq:repfrob}\\
    &=\mathrm{vec}(\mathrm{Var}(Q(X)))' \  (D\otimes D)^{-1} \ \mathrm{vec}(\mathrm{Var}(Q(X))).\label{eq:repkron}
    \end{align}
     We will refer later to the alternative representations~\eqref{eq:repfrob} and \eqref{eq:repkron}. 
\end{Remark}

 The following result confirms that~$\psi(X,Y)$ constitutes a valid measure of dependence that satisfies the key characteristics we demand.
\begin{Theorem}\label{thm:theoretical-properties}
    The dependence coefficient~${\psi(X,Y)}$ satisfies properties~\ref{it:m1}--\ref{it:m3}. Moreover,~$\psi(X,Y)=\psi(X,\ell(Y))$ for any~$\ell\in\mathcal{L}_Y$.
\end{Theorem}

\subsection{Consistent estimation}
Clearly, direct computation of~$\psi(X,Y)$ is practically infeasible even if the joint distribution of~$(X,Y)$ were known. This raises the question of how to estimate~$\psi$ from a random sample~${(X_1, Y_1), \ldots, (X_n, Y_n)}$ drawn from the distribution of~$(X, Y)$. The most natural estimator is obtained by replacing the probabilities~$p_{i}$,~$q_{j}$, and~$p_{i,j}$ with their empirical counterparts. However, since the random variable~$Y'$ is not observable, we need some workaround. The idea is to seek a suitable proxy for~$Y'$. To this end, we use the representation~$Y_i=f(X_i,U_i)$ and~$Y_i'=f(X_i,U_i')$ from~\eqref{eq:y-representation}. We introduce the variable~$N(i)$, which denotes the nearest neighbour of~$X_i$ within the sample of covariates and with respect to the metric~$d$ on~$\mathcal{H}$. 
We set~$U_i'\vcentcolon=U_{N(i)}$ and remark that it is an i.i.d.\ copy of~$U_i$. Moreover, one can show that~$X_{N(i)}\convas X_i$. (See Lemma~\ref{lem:neighbour-convergence} in the supplement.) 
This suggests to replace the unobservable sample~$Y_1',\ldots, Y_n'$ by the proxies~$Y_{N(1)},\ldots, Y_{N(n)}$. 
Let
\begin{equation*}
    \widehat p_{k,l}=\frac{1}{n}\sum_{i=1}^n 1\lbrace Y_i=k,Y_{N(i)}=l\rbrace,\quad  \widehat p_{k} = \sum_{l=1}^K \widehat p_{k,l}\quad\text{and} \quad \widehat q_{l} = \sum_{k=1}^K \widehat p_{k,l}.
\end{equation*}
We then set
\begin{equation*}
    \widehat \psi(X,Y) = \frac{1}{K-1}\sum_{k,l=1}^K \frac{(\widehat p_{k,l}-\widehat p_{k}\widehat q_{l})^2}{\widehat p_{k}\widehat q_{l}}.
\end{equation*}

\begin{Remark}
Nearest neighbours have also been used by   \cite{azadkia_simple_2021} in the construction of their estimator for $\xi$.
\end{Remark}
\begin{Remark}
We have been tacitly assuming that the nearest neighbours are unique. Essentially this means that the distance between two variables $d(X_1,X_2)$ has a continuous distribution.  From a technical point of view uniqueness of $N(i)$ is not important for our results, but it allows to avoid some extra discussions (e.g.\ how to break ties) which are not essential for conveying the main ideas. Hence, for the rest of this paper this assumption is supposed to hold.
\end{Remark}
\begin{Remark}
    The computationally most intensive part in the calculation of~$\widehat\psi$ is finding the nearest neighbour indices~${N(i)}$. Since (approximate) nearest neighbour search can be done in~${O(n\log n)}$ time, the computation time for~$\widehat\psi$ is~${O(n\log n)}$ as well.

\end{Remark}

Our next target is to show consistency of this estimator, i.e.\ convergence to~$\psi(X,Y)$. To this end, we introduce the following high-level assumption.

\begin{Assumption}\label{ass:convergence} Assume that one of the following conditions holds:\\ (a) For almost all~$u$ we have
$
f(X_{N(1)},u)\stackrel{P}{\to} f(X_1,u)$ for $n\to\infty$.\\ (b) For all~$k=1,\ldots,K$, the mapping~$x\mapsto Q_k(x)$ is~$P_X$-almost everywhere continuous.
\end{Assumption}

\begin{Theorem}\label{thm:convergence} Suppose that Assumption~\ref{ass:convergence}  holds. Then~$\widehat\psi(X,Y)$ converges in probability to~${\psi(X,Y)}$.
\end{Theorem}

Given that~$X_{N(1)}$ converges almost surely to~$X_1$,  Assumption~\ref{ass:convergence} seems not to be particularly restrictive. In essence, it allows us to invoke a continuous mapping argument.  We give two sufficient conditions for Assumption~\ref{ass:convergence}. 

\begin{Lemma}
\label{le:suff1}
Assume that~$ \mathcal{H}=\mathbb{R}^d$, then  Assumption~\ref{ass:convergence} holds.
\end{Lemma}

\begin{Lemma}
\label{le:suff2}
Let~$\mathcal{V}_k=\{(x,u)\in \mathcal{H}\times (0,1)\colon f(x,u)=k\}$. If we have \begin{equation}
\label{eq:condsuff2}
P\left((X,U)\in \bigcup_{k=1}^K\partial\mathcal{V}_k\right)=0,
\end{equation}
then Assumption~\ref{ass:convergence} holds.
\end{Lemma}

Condition \eqref{eq:condsuff2} requires that the probability for the pair~$(X,U)$ lying on the discriminant boundaries be zero. This appears to be a natural condition. Since by Lemma~\ref{le:suff1} no extra assumption is needed in finite dimensional scenarios, this condition is only relevant for infinite-dimensional $X$.

In our next result we show that strong consistency can also be achieved. 
Denote with~$\mathcal{G}_n$ the nearest neighbour graph of our covariate sample~$X_1,\ldots, X_n$ and let~$L_n$ be the maximal in-degree of~$\mathcal{G}_n$, i.e.~$L_n=\max_{i=1,\ldots, n} L_{i,n}$, where~$L_{i,n}$ is the number of points that have~$X_i$ as their nearest neighbour.

\begin{Theorem}\label{thm:convergenceas} Let Assumption~\ref{ass:convergence} hold and assume that there is a sequence~$(k_n)$ such that
\begin{equation}\label{ref:assasconv}
\sum_{n\geq 1}P(L_n\geq k_n)<\infty,\quad k_n P(L_n\geq k_n)\to 0 
\quad\text{and}\quad
k_n=o\left(\sqrt{n/\log n}\right),
\end{equation}
        then the estimator~$\widehat\psi(X,Y)$ is strongly consistent. 
\end{Theorem}
\begin{Remark}
    In finite dimensional spaces,~$L_n$ is known to be bounded by a constant which depends only on the dimension of $\mathcal{H}$, and thus Condition~\eqref{ref:assasconv} is trivially satisfied. Along with Lemma~\ref{le:suff1} we conclude that our estimator is strongly consistent if~$\mathrm{dim}(\mathcal{H})<\infty$.
\end{Remark}

\section{Testing for independence}\label{sec:testing}
Associated with Cram\'er's~$V$ is Pearson's~$\chi^2$-test for independence. Suppose for the moment that the~$Y_i'$ were observable and set
$$
\widetilde V(Y,Y')=\frac{1}{K-1}\sum_{k,l=1}^K \frac{(\widetilde p_{k,l}-\widetilde p_{k}\widetilde q_{l})^2}{\widetilde p_{k}\widetilde q_{l}},
$$
with 
\begin{equation*}
    \widetilde p_{k,l}=\frac{1}{n}\sum_{i=1}^n 1\lbrace Y_i=k,Y_i'=l\rbrace,\quad  \widetilde p_{k} = \sum_{l=1}^K \widetilde p_{k,l}\quad\text{and} \quad \widetilde q_{l} = \sum_{k=1}^K \widetilde p_{k,l}.
\end{equation*}
Then it is a well-known result that under independence of~$Y_i$ and~$Y_i'$  it holds that Pearson's~$\chi^2$ statistic
$$
\chi^2\vcentcolon=n(K-1)\widetilde V(Y,Y')
$$
follows an asymptotic~$\chi^2$-distribution with~$(K-1)^2$ degrees of freedom. 

Proceeding similarly as in the previous chapter, we aim to replace the~$Y_i'$ with~$Y_{N(i)}$, or equivalently replacing~$\widetilde{V}(Y,Y')$ with~$\widehat{\psi}(X,Y)$. However, in order to obtain a pivotal limit distribution, we need to adapt this test statistic and account for dependence of~$(Y_i,Y_{N(i)})$ and ~$(Y_j,Y_{N(j)})$ whenever
\begin{equation*}
    \{i,N(i)\}\cap \{j,N(j)\}\neq \emptyset.
\end{equation*}
 Moreover, we must ensure that this dependence effect is not too strong for the envisaged weak convergence. This latter problem can again be controlled by a high-level assumption on~$L_n$:

\begin{Assumption}\label{ass:neighbour-degree}
We have~$L_n=o_P(n^{1/4})$.
\end{Assumption}

\begin{Remark}
    While~$L_n$ remains bounded if~$X$ takes values in~$\mathbb{R}^p$, this is typically not the case in infinite dimensional vector spaces. Indeed, as demonstrated in Theorem~2 of \cite{hormann_azadkiachatterjees_2024}, one can give explicit examples of distributions for which the sequence~$L_n$ diverges at a polynomial rate, but Assumption~\ref{ass:neighbour-degree} is satisfied.
\end{Remark}

\begin{Theorem}\label{thm:clt}
    Assume that~$X$ and~$Y$ are independent, set
\begin{equation*}
P_n=\sqrt{n}\,\mathrm{vec}\Big(((\widehat p_{k,l} -\widehat p_{k}\widehat{ q_{l}}))_{k,l=1}^{K-1}\Big),
\end{equation*}
     and let ~$\Sigma_n=\mathrm{Var}(P_n|X_1,\ldots,X_n)$.
     
\noindent     (i) Then under Assumption~\ref{ass:convergence} 
     \begin{align*}
&\Sigma_n((k_1,l_1), (k_2,l_2))\\
&\qquad=p_{k_1}p_{l_1}p_{k_2}p_{l_2}(1+W_{n})\\
        &\qquad\quad- p_{k_1}p_{l_1}p_{l_2}(1\lbrace k_1=k_2\rbrace + W_{n} 1\lbrace l_1=k_2\rbrace)\\
        &\qquad\quad-p_{k_1}p_{l_1}p_{k_2}(1\lbrace l_1=l_2\rbrace + W_{n} 1\lbrace k_1=l_2\rbrace)\\
        &\qquad\quad+p_{k_1}p_{l_1}(1\lbrace k_1=k_2, l_1=l_2\rbrace + W_{n} 1\lbrace k_1=l_2, l_1=k_2\rbrace)+o_P(1), 
     \end{align*}
     with
\begin{equation}\label{eq:W1n}
W_{n}\vcentcolon=\frac{1}{n}\sum_{i=1}^n I\{N(N(i))=i\}.
\end{equation}
(ii) Let Assumption~\ref{ass:neighbour-degree} hold and suppose that we have some~$b>0$ such that
\begin{equation}\label{eq:lambdamin}
P(W_n\leq 1-b)\to 1.
\end{equation}
Then for large enough~$n$ the matrix~$\Sigma_n$ is invertible and
    \begin{equation*}
    \mathcal{I}_n^0\vcentcolon=P_n'\Sigma_n^{-1}P_n \convd \chi^2((K-1)^2),\quad (n\to\infty).
    \end{equation*}
     If~$\widehat \Sigma_n$ denotes the empirical analogue to~$\Sigma_n$ with~$p_k$ replaced by their estimators~$\widehat p_k$, then we also have 
\begin{equation*}
    \mathcal{I}_n\vcentcolon=P_n'\widehat\Sigma_n^{-1}P_n \convd \chi^2((K-1)^2).
    \end{equation*}
\end{Theorem}

\begin{Remark}
 To clarify the relationship between~$\mathcal{I}_n$ and~$\widehat\psi(X,Y)$, we refer to the representation of~$\psi(X,Y)$ given in \eqref{eq:repkron}. Observe that~$P_n/\sqrt{n}$ serves as the sample analogue of~$\mathrm{vec}(\mathrm{Var}(Q^{[-K]}(X)))$, where~$Q^{[-K]}(X)$ excludes the redundant component
 \begin{equation*}
     Q_K(X) = 1 - Q_1(X) - \cdots - Q_{K-1}(X)
 \end{equation*}
 from~$Q(X)$. To account for dependencies between the coordinates of~$P_n$ and to obtain a pivotal limit distribution, the weighting matrix of the quadratic form is suitably adjusted.
\end{Remark}

The next lemma shows that if the~$X_i$ take values in~$\mathbb{R}^d$, then the matrix~$\Sigma_n$ converges to a non-stochastic limit~$\Sigma$ which is full rank.

\begin{Lemma}\label{lem:limitW1n}
    Suppose~$\mathcal{H}=\mathbb{R}^d$ and assume that~$X_i$ have a continuous density. Then we have~$W_{n}\convP \gamma_d\vcentcolon=\frac{\mathrm{vol}(\mathcal{S}_d)}{\mathrm{vol}(\mathcal{S}_d\cup 
    \mathcal{S}_{1,d})}$, where~$\mathcal{S}_d$ is the unit-sphere in~$\mathbb{R}^d$ and~$\mathcal{S}_{1,d}$ is an affine shift of~$\mathcal{S}_d$ by distance 1. It holds that~$\gamma_d\in (1/2, 2/3]$, with~$\gamma_1=2/3$ and~$\gamma_d\searrow  1/2$ for~$d\to\infty$.
\end{Lemma}
Lemma~\ref{lem:limitW1n} follows from general results of \cite{henze_fraction_1987} (Remark~1.1) and \cite{schilling_mutual_1986} (Lemma~4.1) about nearest-neighbour interrelations. The latter paper also contains a table for~$\gamma_d$, with selected  values of $d$.

The following corollary highlights an important special case of Theorem~\ref{thm:clt}.

\begin{Corollary}
Suppose that~$Y$ is a binary variable and that~$X$ takes values in~$\mathbb{R}^d$. Then if~$X$ and~$Y$ are independent we have
$$
\mathcal{I}_n^d\vcentcolon=\frac{n}{1+\gamma_d}\left(\frac{\widehat p_{1,1}-\widehat p_1\widehat q_1}{\widehat{p}_1\widehat{q}_1}\right)^2\convd \chi^2(1).
$$
\end{Corollary}

 If the target is to test
\begin{equation}
\text{
$\mathcal{H}_0\colon X_i$ and~$Y_i$ are independent\quad v.s.\quad~$\mathcal{H}_A\colon \mathcal{H}_0$ doesn't hold,
}
\end{equation}
then the test rejecting~$\mathcal{H}_0$ if~$\mathcal{I}_n>q_{1-\alpha}\big(\chi^2_{(K-1)^2}\big)$ has asymptotic significance level~$\alpha$ and can be performed in~${O(n\log n)}$ time, without any resampling. The following corollary shows that under our assumptions, this test is universally consistent.

\begin{Corollary}\label{cor:consistent-test}
Let Assumptions~\ref{ass:convergence} and~\ref{ass:neighbour-degree} hold. Assume that~$X$ and~$Y$ are not independent. Then
$\mathcal{I}_n\convP\infty$  as~$n\to\infty$. The divergence rate is of order  $\sqrt{n}$.
\end{Corollary}

\begin{Remark}
    As outlined in Section~\ref{sec:introduction}, testing independence between $X$ and $Y$ is equivalent to testing whether the $K$ subsamples of $X$ associated with different levels of $Y$ originate from the same distribution. The test proposed above therefore also contributes to the $K$-sample problem, with particular relevance in the challenging setting of high- or even infinite-dimensional distributions.
\end{Remark}

\section{Further properties and extensions}\label{sec:extensions}
A dependence measure satisfying~\ref{it:m1}-\ref{it:m3} can be seen as a measure of \textit{predictability}. As noted by~\cite{ansari_direct_2025}, such a measure of predictability should have the property that additional information increases the predictability of~$Y$. This requirement is made formal by
\begin{enumerate}[label = (M\arabic*), start=4]
\item\label{it:m4}~$\psi(X,Y)\leq \psi((X,Z),Y)$ for all~${X,Z,Y}$ (\textit{information gain inequality}) and
    \item\label{it:m5}~$\psi(X,Y)= \psi((X,Z),Y)$ if and only if~$Y$ and~$Z$ are conditionally independent given~$X$.
\end{enumerate}
The following result, which will be important for the subsequent discussion, confirms that~$\psi$ satisfies these properties. 

\begin{Theorem}\label{thm:information-gain}
    The dependence coefficient~${\psi(X,Y)}$ satisfies~\ref{it:m4}--\ref{it:m5}.
\end{Theorem}

\subsection{Measuring conditional dependence}\label{sec:conditional}
We propose an extension of the measure~$\psi$ to a measure of conditional dependence: Let~$Z$ be a random variable with values in a separable metric space~$\mathcal H'$. We define
\begin{equation}\label{eq:conditional}
    \psi(Z,Y|X)\vcentcolon=\frac{\psi((X,Z),Y)-\psi(X,Y)}{1-\psi(X,Y)},
\end{equation}
if~$\psi(X,Y)\neq 1$. If~$\psi(X,Y)= 1$, a measure of conditional dependence is pointless, as all information about~$Y$ is already contained in~$X$ and therefore~$Z$ cannot provide any additional information. In analogy to~\ref{it:m1}--\ref{it:m3}, the following properties are desirable:
\begin{enumerate}[label = (M\arabic*')]
        \item\label{it:m1-conditional}~$0\leq \psi(Z,Y|X)\leq 1$,
        \item\label{it:m2-conditional}~$\psi(Z,Y|X)=0$ if and only if~$Z$ and~$Y$ are independent conditional on~$X$,
        \item\label{it:m3-conditional}~$\psi(Z,Y|X)=1$ if and only if there is a measurable function~$f$ such that~${Y=f(X,Z)}$.
    \end{enumerate}

\begin{Theorem}\label{thm:conditional}
    The conditional dependence measure~$\psi(Z,Y|X)$ satisfies properties~\ref{it:m1-conditional}--\ref{it:m3-conditional}. Moreover,~${\psi(Z,Y|X)=\psi(Z,\ell(Y)|X)}$ for any~${\ell\in\mathcal{L}_Y}$. 
\end{Theorem}
The first statement of Theorem~\ref{thm:conditional} is a special case of the more general Theorem~\ref{thm:general-norms} below. The permutation invariance is directly inherited from the unconditional variant.
\begin{Remark}
    If the conditions of Theorems~\ref{thm:convergence} or~\ref{thm:convergenceas} are fulfilled (for~${\mathcal H, \mathcal{H}'}$), the conditional dependence measure~${\psi(Z,Y|X)}$ can be consistently estimated by plugging the estimators for~${\psi((X,Z),Y)}$ and~${\psi(X,Y)}$ into~\eqref{eq:conditional}.
\end{Remark}
A measure of conditional dependence satisfying~{(M2')} can be naturally used for the task of variable selection, since it does not need any specific model assumption (see Section~5 of \cite{azadkia_simple_2021} for a more detailed discussion): Given a set of possible covariates~${X_1,\ldots,X_p}$ for a response~$Y$, one selects the first variable as
\begin{equation*}
    i_1 =\argmax_{1\leq i\leq p} \psi(X_i,Y).
\end{equation*}
Inductively, one selects the~$j$-th variable as
\begin{equation*}
    i_j =\argmax_{\substack{1\leq i\leq p\\ i\neq i_1,\ldots,i_{j-1}}} \psi(X_i,Y|X_{i_1},\ldots,X_{i_{j-1}}),
\end{equation*}
until some stopping criterion is fulfilled. In our simulations we use the smallest~$j$ for which~$\widehat{\psi}(X_{i_j},Y|X_{i_1},\ldots, X_{i_{j-1}})\leq 0$.

\begin{Remark}
One can also use $\widehat\psi$ to get a test for conditional dependence by using the conditional randomization testing \citep{candes_panning_2018} or conditional permutation testing \citep{berrett_conditional_2020} frameworks.
\end{Remark}

\subsection{Alternative norms}
As stated in \eqref{eq:repfrob}, the dependence measure~$\psi$ can be expressed as the squared (weighted) Frobenius norm of~$\Var(Q(X))$. This observation can be used to define a general class of dependence measures. To this end, we recall that a matrix norm~${\|.\|}$ is \textit{monotone}, if for all symmetric, positive semidefinite (PSD), real matrices~${A,B}$, we have
    \begin{equation*}
        B-A~\text{is PSD} \Rightarrow \|A\|\leq \|B\|.
    \end{equation*}
We call~${\|.\|}$ \textit{strictly monotone}, if
    \begin{equation*}
        B-A~\text{is PSD and}~B\neq A\Rightarrow \|A\|<\|B\|.
    \end{equation*}
Examples of strictly monotone matrix norms include the Frobenius norm and the trace norm. The spectral norm is an example of a monotone matrix norm that is not strictly monotone.

\par Given a strictly monotone matrix norm~${\|.\|}$ and a strictly monotone function~${\gamma\colon\lbrack 0,1\rbrack\to\lbrack0,1\rbrack}$ satisfying~${\gamma(0)=0}$ and~${\gamma(1)=1}$, we define 
\begin{equation}\label{eq:norm-measure}
    \psi_{\|.\|,\gamma}(X,Y)\vcentcolon=\gamma\left(\frac{\|\Var(Q(X))\|}{\|\Var(Q(Y))\|}\right).
\end{equation}
The factor~${\|\Var(Q(Y))\|^{-1}}$ serves to normalize the coefficient to~1 for the case of perfect functional dependence (since~$Y$ is obviously a function of~$Y$). The measure~$\psi$ arises for the choice of the weighted Frobenius norm $$\|A \|_{F,D}=\|D^{-1/2}A D^{-1/2}\|_{F}$$ and~${\gamma\colon x\mapsto x^2}$. In this case, it is easily seen that~${\|\Var(Q(Y))\|^2_{F,D}=K-1}$. Similarly, it can be seen that for the weighted trace norm
\begin{equation*}
    \| A\|_{\mathrm{tr}, D} = \mathrm{trace}(D^{-1/2}AD^{-1/2}),
\end{equation*}
we have~${\|\Var(Q(Y))\|_{\mathrm{tr}, D} = K-1}$. In case~${\|\Var(Q(Y))\|}$ is unknown for a given norm, it can be easily estimated, since~${\Var(Q(Y))}$ has entries
\begin{align}\label{eq:var-q}
    a_{kl} =\begin{cases}
        p_k-p_k^2&\text{if}~k=l\\
        -p_kp_l&\text{otherwise}.
    \end{cases}
\end{align}
We may then simply plug in the estimators~$\widehat p_k$ into~\eqref{eq:var-q} and calculate the norm of the resulting matrix to obtain a strongly consistent estimator.
\par Obviously, any such measure can be used to measure conditional dependence as well by using~\eqref{eq:conditional}.
\par We conclude by showing that every measure $\psi_{\|.\|,\gamma}(X,Y)$ has the same theoretical properties as~$\psi(X,Y)$.
\begin{Theorem}\label{thm:general-norms}
    If~${\|.\|}$ is a strictly monotone norm and~$\gamma$ is a function as described above, then the measure~$\psi_{\|.\|,\gamma}$ satisfies the properties~\ref{it:m1}--\ref{it:m5} as well as~\ref{it:m1-conditional}--\ref{it:m3-conditional}.
\end{Theorem}

\begin{Remark}
    Theorem~\ref{thm:general-norms}  provides  some intuition about  the nature of our approach. Determining whether alternative norms yield more powerful independence tests necessitates an investigation into the large-sample behaviour of these coefficients. This analysis lies beyond the scope of the current study and is left for future research.
\end{Remark}

\section{Data illustration}\label{sec:illustration}
In this section, we evaluate the performance of the proposed dependence measure and independence test on simulated and real-world data. Implementations of the calculation of~$\widehat\psi$, the independence test, and the variable selection algorithm are available through the package {\tt R}-package \texttt{FDEP} from from the second author's {\tt GitHub} profile \citep{strenger-galvis_fdep_2025}.
\subsection{Simulations}

We compare the power of our test ({\tt psicor}, $\widehat\psi$) in various settings to the power of tests based on distance correlation ({\tt dcor}), ball correlation ({\tt bcor}), Gini distance correlation ({\tt gcor}) and Generalized Maximum Mean Discrepancy ({\tt gmmd}). For details on these competing tests, we refer to \cite{szekely_distance_2013}, \cite{pan_ball_2020}, \cite{dang_new_2021}, and \cite{balogoun_k-sample_2022}, respectively. We consider different types of noisy functional relationships between~$X$ and~$Y$. 
Specifically, we consider
\begin{itemize}
    \item \quad $(X_1,X_2)\sim\mathcal{U}(\lbrack0,1\rbrack^2)$ and~$Y=I\lbrace\sin(2\pi(X_1+X_2))\geq0\rbrace$. ({\tt sin})
\item \quad $X$ is a standard Brownian motion and~${Y=I\lbrace \max_{t\in\lbrack0,1\rbrack}X(t)\geq 1\rbrace}$. ({\tt max})
    \item \quad  The variable~$X$ is drawn from two functional populations: either~$X$ is a Brownian motion or it is given as~${X(t)=\sum_{k=1}^{20}} Z_k\sin(\pi kt)$, with random variables~${Z_k\sim\mathcal N(0,0.5^k)}$. The respective population is expressed by the categorical variable~$Y$. Each of the two populations has prior probability~$1/2$. ({\tt mixture})
    
    \item \quad $X$ is a random polynomial of degree~$Y$, where~$Y$ is drawn uniformly from~${0,\ldots,8}$. The nonzero coefficients of the polynomial are drawn uniformly and independently from the interval~${\lbrack0,1\rbrack}$. ({\tt degree})
\end{itemize}

\begin{figure}[h]
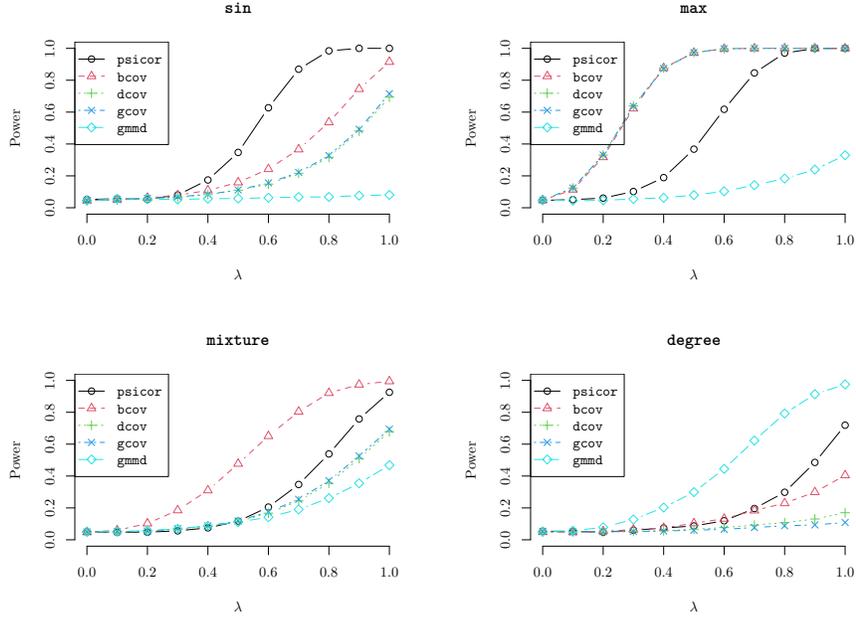

       \resizebox{0.49\textwidth}{!}{
        \input{powers-sin-gmmd}
        }
    \resizebox{0.49\textwidth}{!}{
        \input{powers-max-gmmd}
        }\\     
    \resizebox{0.49\textwidth}{!}{
        \input{powers-mixture-gmmd}
        }
     \resizebox{0.49\textwidth}{!}{
        \input{powers-degree-gmmd}
        }
    \caption{Power functions of the different independence tests. The case $\lambda=0$ corresponds to the null-hypothesis of independence between $Y$ and $X$, while~${\lambda=1}$ corresponds to a functional relation $Y=f(X)$.}
    \label{fig:powers}
\end{figure}

In all these examples~$Y$ is measurable with respect to~$X$.  In order to induce nondeterministic relationships, we calculate the test statistics for the pairs~${( X,\widetilde Y)}$, where
\begin{equation*}
    \widetilde Y=\delta Y+(1-\delta) Y',
\end{equation*}
with~$Y'$ being an i.i.d.\ copy of~$Y$ and~$\delta$ being a binomial~$B_{1,\lambda}$ variable, independent of everything else. Testing for independence of $Y$ and $X$ is equivalent to testing~${\mathcal{H}_0\colon \lambda=0}$. On the other end, when~$\lambda=1$, we have complete dependence. In each setting, we draw~10,000 samples of size~$n=100$ to estimate the power as function in~$\lambda$. The number of repetitions for the resampling based {\tt dcor}, {\tt bcor}, {\tt gcor} and {\tt gmmd} tests was set to~${R=500}$.  The results are given in Figure~\ref{fig:powers}.

All tests attain the correct size (the nominal level set to~$5\%$). Within the five considered setups there is no method dominating the others. Our test performs best in settings {\tt sin}. The methods {\tt dcov} and {\tt gcov} yield almost identical results.

As discussed above, the test based on \texttt{psicor} relies on the asymptotic distribution instead of a resampling procedure, which significantly reduces the computation time. In addition, a single computation of the test statistic takes~${O(n\log n)}$ time, compared to~$n^2$ for \texttt{dcor} \texttt{gcor} and \texttt{gmmd}, and~$n^3$ for \texttt{bcor}. Table~\ref{tab:running-times} compares the running times for a bivariate covariate $X$ and a response~$Y$ with three levels. The number of repetitions for the resampling-based tests was set to~${R=500}$. The computation times were measured on a desktop computer containing an~i5-11600 CPU.  For larger datasets, the competing methods require extensive computation time.  This makes the \texttt{psicor} based test well suited for situations in which large amounts of data are available. With sufficiently large sample size, the power of the test approaches~1 (see Corollary~\ref{cor:consistent-test}). The competing methods might be inapplicable in these situations due to their computational cost.

\begin{table}[]
    \centering
    \caption{Running times (in seconds) of the independence tests for different sample sizes.}
    \label{tab:running-times}
    \begin{tabular}{rrrrr}
    \toprule
         $n$& \texttt{psicor} & \texttt{dcor} & \texttt{bcor} & \texttt{gcor}  \\
         \midrule
         100 & $<0.01$ & $<0.01$ & $0.22$ &  $0.19$\\
         500 & $<0.01$ & $0.12$ & $5.08$ & $1.95$\\
         1000 & $<0.01$ & $0.66$ & $21.25$ & $6.90$\\
         2000 & $<0.01$ & 3.14 & 127.94 & 38.42\\
         10000 & 0.01 & 114.30 & 3354.94 & 1060.47
    \end{tabular}
\end{table}

\subsection{Real data}
In the following we consider different real data illustrations.
\par
\subsubsection{Election data} 
For each of the~${n=2061}$ Austrian municipalities, we have a vector~$X$ containing some demographic and socio-economic indicators.  The variables include age distribution (percentage of the population within predefined age groups: 19–25, 30–44,~\ldots,~${\geq75}$ years), educational attainment (percentage of the population by highest completed qualification, ranging from mandatory schooling to tertiary education), the share of non-EU citizens in the population, and population density (inhabitants per square kilometer). This data is provided by \cite{statistik_austria_statatlas_2025}.
 The categorical variable~$Y$, provided by \cite{sora_national_2019}, indicates which of two candidates of the 2016 presidential election received the most votes in a municipality. Since there were a few smaller municipalities with a tie, the variable~$Y$ is defined with~3 levels. The target it explore the strenght of dependence between the socio-economic variables and the selected candidate.
\par The measures of dependence take the values~${\widehat \psi=0.22, \widehat R =0.26, \widehat B=0.03,}$ and~${\widehat G = 0.035}$, indicating that~$X$ has a significant (the~$p$-values of the corresponding independence tests were~$<0.01$) impact on~$Y$, but is far from determining~$Y$ completely.

\subsubsection{Spambase and bankruptcy data} 
\cite{azadkia_simple_2021}
have applied their feature selection procedure FOCI to the \textit{spambase} and \textit{polish company bankruptcy} datasets from the UCI machine learning repository \citep{kelly_uci_2025}. Both have a large number of features (57 and 64, respectively) and a binary response variable (spam/no spam and bankruptcy/no bankruptcy). We do similar  comparisons as \cite{azadkia_simple_2021} with the LASSO and the Dantzig selector \citep{candes_dantzig_2007}. Instead of  step-wise  regression we include Gini Distance Correlation \citep{dang_new_2021}. We choose a random subset comprising 90\% of the data to perform the feature selection procedures and fit a random forest model using the variables selected by each procedure. In Table~\ref{tab:feature-selection} we report the number of selected variables and the prediction accuracy on the remaining~10\% test dataset. 

We note that for our method and for FOCI the variable selection is sequential and based on a stopping criterion: both methods stop, when the respective dependence coefficient between $Y$ and any newly selected variable indicates independence, conditional on the  variables that have been already selected. The other methods were tuned using~10-fold cross-validation. This implies that the variables are specifically selected for achieving high prediction accuracy, which is of course beneficial in this competition. However, this advantage comes at the price of significantly increased computation time for choosing substantially larger models, while the increase in accuracy is relatively small in comparison.

\begin{table}[h]
	\renewcommand\arraystretch{1.2}
	\begin{center}
	\caption{Comparison of different feature selection methods. \label{tab:feature-selection}}
			\begin{tabular}{lcccccc}
				\toprule
				& \multicolumn{2}{c}{Spambase data} & & \multicolumn{2}{c}{Polish companies data}\\
				\cmidrule{2-3} \cmidrule{5-6}
				Method & Subset size& Accuracy & & Subset size & Accuracy \\
			\midrule
			PSICOR & 14 & 0.918 & & 1 & 0.969 \\
			FOCI & 19 & 0.939 & & 6 & 0.981\\
			Lasso & 53 & 0.957  & & 6 & 0.976 \\
            Gini & 57 & 0.959 & & 57 & 0.980\\
			Dantzig selector & 57 & 0.959 & & 29 & 0.976 \\
			Full dataset & 57 & 0.959 & & 64 & 0.979 \\
			\bottomrule
			\end{tabular}
	\end{center}
\end{table}

\par
\subsubsection{MNIST data} 
\begin{figure}[h]
    \centering
    \includegraphics[width=0.5\textwidth]{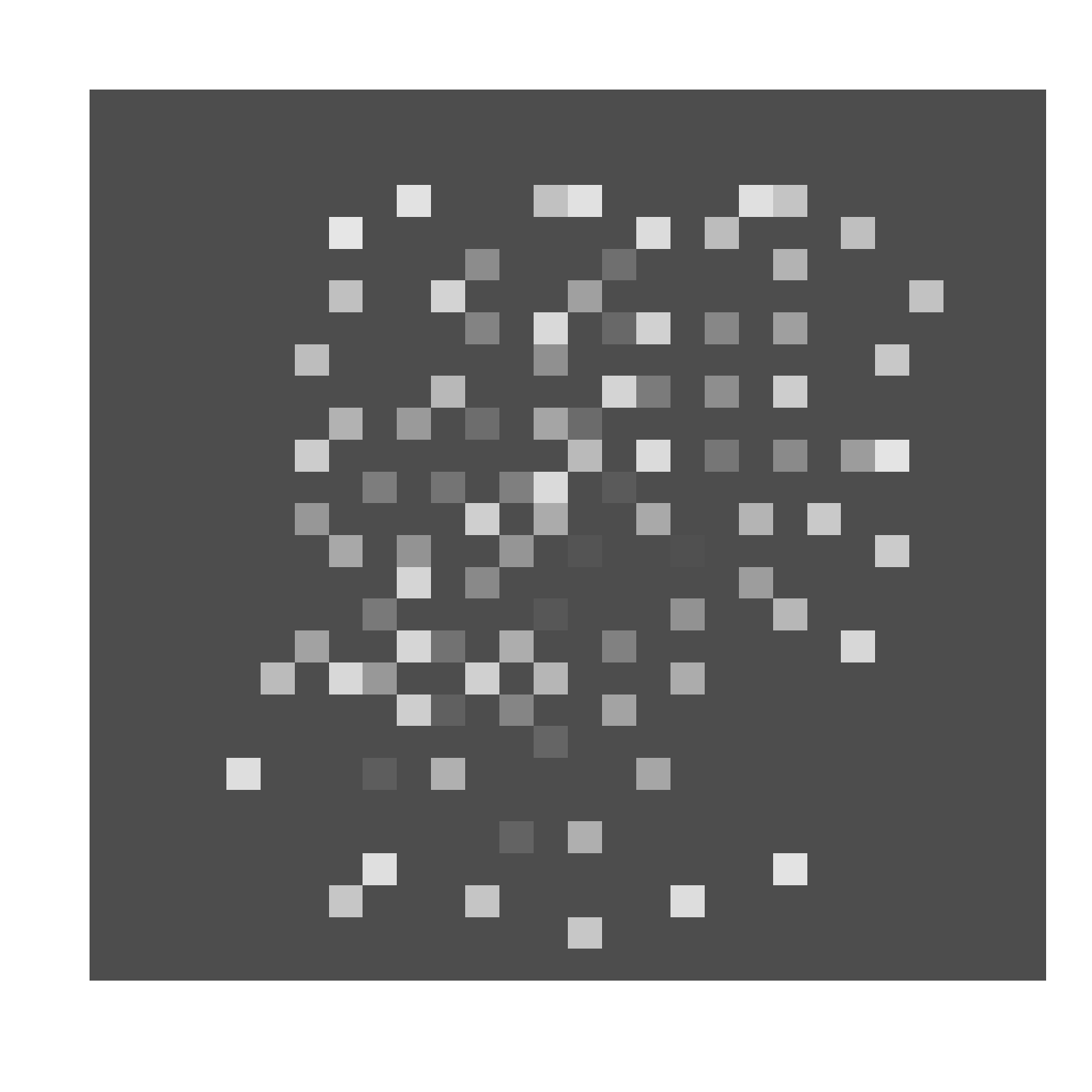}
    \caption{Illustration of the pixels selected by {\tt psicor} feature selection. The brighter the pixels the earlier they have been selected.}
    \label{fig:selected}
\end{figure}

\begin{figure}[h]
    \centering
    \includegraphics[width=\textwidth]{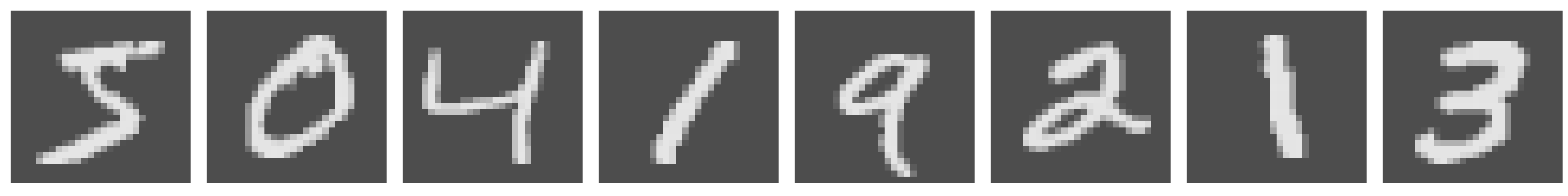}\\
    \includegraphics[width=\textwidth]{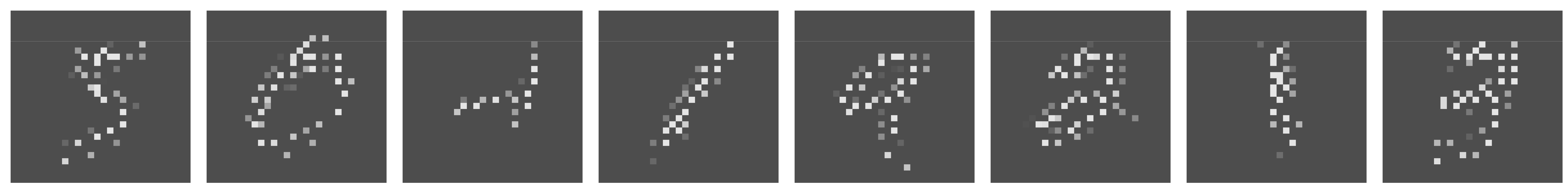}
    \caption{The first eight images from the MNIST training data (top) and their selected pixels (bottom).}
    \label{fig:digits}
\end{figure}

The MNIST dataset \citep{lecun_gradient-based_1998} consists of~${28\times28}$-pixel grayscale images~$X$ of handwritten digits, as well as labels~$Y$ indicating which digit is depicted. As we can visually classify the depicted digit in most cases, it is clear that in theory~${\psi(X,Y)}$ should be close to~1.  This is reflected by the value~${\widehat \psi (X,Y)=0.951}$ on the~10,000 test images. Moreover, one can expect that spatially close pixels are highly dependent and therefore a subset of the~784 pixels should suffice to recognize a digit. This is confirmed by our feature selection algorithm identifying only~96 pixels as relevant. An illustration of the selected pixels can be seen in Figures~\ref{fig:selected} and~\ref{fig:digits}. Most of the digits are clearly recognizable for humans based on only the selected pixels. An XGBoost model \citep{chen_xgboost_2016} trained to identify the digit based on only the selected pixels achieves an accuracy of~0.979, almost the same as the~0.980 achieved by using all pixels.

\section{Proofs}\label{sec:proofs}
The proofs of the auxiliary results presented in this section can be found in Supplements~\ref{sec:appendix} and~\ref{sec:calculation-sigma}.
\subsection{Proof of Proposition~\ref{pr:1}}
We construct an example with a three levels categorical variable~$Y$ and two-levels categorical variable~$X$. For~$\varepsilon\in(0,1)$, consider the joint distribution as given in Table~\ref{tab:joint-dist-gen}.

\begin{table}[h!]
    \centering
    \caption{Joint probability function of~$Y$ and~$X$.}
    \label{tab:joint-dist-gen}
    \renewcommand{\arraystretch}{1.3}
    \begin{tabular}{cccc}
        \toprule
        \backslashbox{$X$}{$Y$} &~$A$ &~$B$ &~$C$ \\
        \midrule
       ~$0$ &~$0$ &~$1-\varepsilon$ &~$\frac{\varepsilon}{2}$ \\
       ~$1$ &~$\frac{\varepsilon}{2}$ &~$0$ &~$0$ \\
        \bottomrule
    \end{tabular}
\end{table}

Using the encoding~$A=1, B=2, C=3$, some basic computations show
    \begin{align*}
        \int \Var\left( 1\lbrace Y\geq t\rbrace\right) d\mu(t)
        &=\varepsilon/2-\varepsilon^2/2+\varepsilon^3/8
    \end{align*}
    and
    \begin{align*}
        \int\Var\left(P(Y\geq t|X)\right) d\mu(t)
        &= \varepsilon/2-3\varepsilon^2/4+\varepsilon^3/4+\frac{\varepsilon^4}{16(1-\varepsilon/2)}.
    \end{align*}
    Therefore, we get that~$\xi(X,Y)\to 1$, as~$\varepsilon\to 0$.
    \par On the other hand, if we use the encoding~$A=3, B=2, C=1$, 
    we obtain
    \begin{align*}
         \int\Var\left(P(Y\geq t|X)\right) d\mu(t)
         &=\varepsilon^2/4\left(\left(1-\varepsilon/2\right)+\frac{\varepsilon}{2}\frac{1-\varepsilon}{1-\varepsilon/2}\right).
    \end{align*}
    Therefore, we obtain~$\xi(X,Y)\to 0$, as~$\varepsilon\to 0$. 

\begin{Remark}
    Note that in the extreme case~$\varepsilon=0$ we have~$Y=B$ almost surely, and hence a degenerate case. As noted by \cite{chatterjee_new_2021}, in this case~$Y$ is both, independent of~$X$ and equal to a function of~$X$ and the measure~$\xi$ is undefined, as the denominator in~\eqref{eq:xi-definition} is equal to zero. In this specific example, each of the two encodings captures the ``closeness'' to this case of combined full dependence \textit{and} independence from a different perspective. In turn, one can easily derive that~${\psi(X,Y)=1/2}$ for all~${\varepsilon\in (0,1)}$.
\end{Remark}

\subsection{Proof of Theorem~\ref{thm:theoretical-properties}}
\begin{proof}[Proof of Lemma~\ref{lem:x-ind-y}]
    (i) Obviously,  independence of~$X$ and~$Y$ implies independence of~$Y$ and~$Y'$. For the converse implication, consider
    \begin{align}
&P(Y=k, Y'=k)=\int_\mathcal{H} P(Y=k, Y'=k|X=x)dP_X(x)\nonumber\\
&\qquad=\int_\mathcal{H} P(Y=k|X=x)^2dP_X(x)\nonumber\\
&\qquad\geq\left(\int_\mathcal{H} P(Y=k|X=x)dP_X(x)\right)^2\label{e:ineq}\\
&\qquad=\left(\int_\mathcal{H} P(Y=k|X=x)dP_X(x)\right)\left(\int_\mathcal{H} P(Y'=k|X=x)dP_X(x)\right)\nonumber\\
&\qquad=P(Y=k)P(Y'=k).\nonumber
\end{align}
The~$\geq$ in~\eqref{e:ineq} is due to Jensen's inequality, which is strict, unless the integrand is a.s.\ constant. Thus, if~$Y$ and~$Y'$ are independent, then the integrand~${P(Y=k|X=x)}$ must not depend on~$x$, and therefore~$Y$ and~$X$ are independent.

\par
(ii) If~$Y$ is a function of~$X$, then~${Y=f(X)=Y'}$, so the function linking~$Y$ and~$Y'$ is the identity. To show the reverse implication we note that~${Y=\ell(Y')}$ implies that we have~${f(X,U)=\ell(f(X,U'))}$. Since~$X$,~$U$ and~$U'$ are independent, this equality can only hold if~$f$ doesn't depend on the argument~$U$.
\end{proof}

\begin{proof}[Proof of Theorem~\ref{thm:theoretical-properties}]
The invariance of~$\psi$ to permuations of the levels of~$Y$ is immediate from its definition. Property~\ref{it:m1}-\ref{it:m3} are known to hold for Cram\'er's~$V$. The result is then immediate from Lemma~\ref{lem:x-ind-y}.
\end{proof}

\subsection{Proofs of Section~\ref{sec:measure}}
We first state some preparatory lemmas.

\begin{Lemma}\label{lem:convergence}
    (i) Under Assumption~\ref{ass:convergence}  it holds for all~${k,l}$ that
    \begin{equation*}
        \widehat p_{k,l}\convP p_{k,l},\quad (n\to\infty).
    \end{equation*}
    \par (ii) If we additionally impose \eqref{ref:assasconv}, then convergence in probability can be strengthened to almost sure convergence.
\end{Lemma}

\begin{proof}[Proof of Theorems~\ref{thm:convergence} and~\ref{thm:convergenceas}]
    Both theorems follow directly from Lemma~\ref{lem:convergence} and the continuous mapping theorem.
\end{proof}

\begin{proof}[Proof of Lemma~\ref{le:suff1} and Lemma~\ref{le:suff2}]  One can show that for any measurable function~${g\colon \mathbb{R}^p\to\mathbb{R}}$ it holds that
${g(X_{N(1)})\stackrel{P}{\to} g(X_1)}$ 
\cite[Lemma~11.7]{azadkia_simple_2021}. Hence, for multivariate~$X_1$ the condition \eqref{eq:altcond}---which is equivalent to Assumption~{\ref{ass:convergence}~(a)}---holds. By the continuous mapping theorem,~\eqref{eq:condsuff2} also implies~\eqref{eq:altcond}.
\end{proof}


\subsection{Proofs of Section~\ref{sec:testing}}
For the remainder of this section, assume that~$X$ and~$Y$ are independent. Note that this implies~${Y_{N(1)} \stackrel{\mathcal{D}}{=} Y_1}$. We denote by~$\mathcal{G}_n$ the nearest-neighbour graph of~$X_1,\ldots, X_n$.

 We write 
    \begin{align*}
        \sqrt{n}(\widehat p_{k,l}-\widehat p_k\widehat q_l) &=\frac{1}{\sqrt{n}} \sum_{i=1}^n (1\lbrace Y_i =k\rbrace -\widehat p_k)(1\lbrace Y_{N(i)} =l\rbrace -\widehat q_l)\nonumber\\
        &= \frac{1}{\sqrt{n}} \sum_{i=1}^n (1\lbrace Y_i =k\rbrace - p_k)(1\lbrace Y_{N(i)} =l\rbrace - p_l)\\
        &\quad - \sqrt{n}(p_k-\widehat p_k)(p_l-\widehat q_l).
    \end{align*}
    Using the Cauchy-Schwarz inequality and Lemma~\ref{lem:convergence} it can be readily verified that
    \begin{equation*}
        n E(p_k-\widehat p_k)^2(p_l-\widehat q_l)^2\to 0.
    \end{equation*}
    Therefore, for the proof of Theorem~\ref{thm:clt} we may slightly modify the definition of~$P_n$. Redefining 
    \begin{equation*}
P_n\vcentcolon=\mathrm{vec}\left(\left(\left(\frac{1}{\sqrt{n}} \sum_{i=1}^n (1\lbrace Y_i =k\rbrace - p_k)(1\lbrace Y_{N(i)} =l\rbrace - p_l)\right)\right)_{k,l=1}^{K-1}\right)
    \end{equation*}
    induces an asymptotically negligible perturbation with respect to both, the (conditional) covariance matrix and the limiting distribution of~$P_n$.

\begin{proof}[Proof of Theorem~\ref{thm:clt} (i)]  
    Since
    \begin{equation*}
        \ev\left(1\lbrace Y_i =k\rbrace\right)=p_k\qquad\text{and}\qquad \ev\left(1\lbrace Y_{N(i)} =l\rbrace\right)=p_l,
    \end{equation*}
    and~${Y_i, Y_{N(i)}}$ are independent (unconditionally and conditionally on~$\mathcal{G}_n$), we get that
    \begin{equation*}
        \ev\left((1\lbrace Y_i =k\rbrace - p_k)(1\lbrace Y_{N(i)} =l\rbrace - p_l)|\mathcal{G}_n\right) = 0.
    \end{equation*}
    As a consequence, it follows that
    \begin{align*}
    &\Cov\Big((1\lbrace Y_i=k_1\rbrace -p_{k_1})(1\lbrace Y_{N(i)}=l_1\rbrace -p_{l_1}),\\ &\qquad\qquad (1\lbrace Y_j=k_2\rbrace -p_{k_2})(1\lbrace Y_{N(j)}=l_2\rbrace -p_{l_2})|\mathcal{G}_n\Big)\\
        &= \ev\big((1\lbrace Y_i=k_1\rbrace -p_{k_1})(1\lbrace Y_{N(i)}=l_1\rbrace -p_{l_1})\times\\
        &\qquad\qquad
        (1\lbrace Y_j=k_2\rbrace -p_{k_2})(1\lbrace Y_{N(j)}=l_2\rbrace -p_{l_2})|\mathcal{G}_n\big).
    \end{align*}
    Expanding these products, we obtain~16 summands, depending the combinations of~$i$ and~$j$. We defer these tedious calculations to Supplement~\ref{sec:calculation-sigma}. Summing over these terms yields~$\Sigma_n$.
\end{proof}

We proceed to prove the CLT for our test statistic. We will make use of the following normal approximation theorem.

\begin{Lemma}\label{lem:clt}
    We denote~$\Sigma_{g,n}=\mathrm{Var}(P_n|\mathcal{G}_n=g)$ and let~$\mathcal{K}_n, n\geq 1$, be a sequence of collections of nearest neighbour graphs~$g$ with~$n$ vertices, such that
    \begin{equation*}
        \sup_{g\in \mathcal{K}_n}\mathrm{deg}(g)= o\left(n^{1/4}\right)
    \end{equation*}
    and such that the coefficient~$W_{n}$ defined in \eqref{eq:W1n} and corresponding to~$g$ is at most~${1-b/2}$. Then, for any vector~$c$ on the~$(K-1)^2$-dimensional unit sphere it holds that
    \begin{equation*}
        \sup_{g\in\mathcal K_n}\sup_{z\in\R} \left|P\left(c'\Sigma_{g,n}^{-1/2}P_n\leq z\middle| \mathcal G_n = g\right)-\Phi(z)\right|\to 0,\quad (n\to \infty).
    \end{equation*}
\end{Lemma}

\begin{proof}[Proof of Theorem~\ref{thm:clt}, (ii)]
Assumptions~\ref{ass:neighbour-degree} and \eqref{eq:lambdamin} imply that~$P(\mathcal{G}_n\in\mathcal{K}_n)\to 1$. Hence, by the Wold device and Lemma~\ref{lem:clt} we have
that~$\Sigma_n^{-1/2}P_n$ converges in distribution to~$(K-1)^2$-variate standard normal vector. The asymptotic~$\chi^2$-distribution of our test statistic follows by the continuous mapping theorem.

Finally, condition~\ref{eq:lambdamin} and Lemma~\ref{lem:convergence} assure that~$\widehat\Sigma_n^{1/2}\Sigma_n^{-1/2}$ converges in probability to the identity matrix, which finishes our proof.
\end{proof}

\begin{proof}[Proof of Corollary~\ref{cor:consistent-test}]
By Lemma~\ref{lem:x-ind-y}, there are indices~${k,l}$ with~${p_{k,l}\neq p_kp_l}$. By Lemma~\ref{lem:convergence}, it follows that
    \begin{equation*}
        |\sqrt{n} (\widehat p_{k,l}-\widehat p_k\widehat q_l)|\convP \infty.
    \end{equation*}
    Since~$\widehat\Sigma_n$ is asymptotically regular, the test statistic diverges as well.
\end{proof}

\subsection{Proof of Theorem~\ref{thm:general-norms}}
\begin{Lemma}\label{lem:conditional}
Under the assumptions of Theorem~\ref{thm:general-norms} it holds that
\begin{equation*}
    \|\Var(Q(X,Z))\|\geq \|\Var(Q(X))\|,
\end{equation*}
with equality if and only if~$Y$ is independent of $Z$ conditionally on $X$.
\end{Lemma}
\begin{proof}
Since~${Q(X,Z)=(P(Y=1|X,Z),\ldots, P(Y=K|X,Z))'}$, it holds that
\begin{equation*}
    \ev(Q(X,Z)|X)=Q(X).
\end{equation*}
By the law of total variance, we have
\begin{align*}
    &\|\Var(Q(X,Z))\|\\
    &\quad = \|\Var(\ev(Q(X,Z)|X)) + \ev(\Var(Q(X,Z)|X))\|\\
    &\quad = \|\Var(Q(X)) + \ev(\Var(Q(X,Z)|X))\|\\
    &\quad\geq \|\Var(Q(X))\|,
\end{align*}
because~${\|\cdot\|}$ is monotone and (conditional) covariance matrices as well as their expectations are PSD. Moreover, by strict monotonicity of~${\|\cdot\|}$, equality holds if and only if~${\Var(Q(X,Z)|X)}$ is almost surely equal to the zero matrix. But this is only the case if~$Q(X,Z)$ is~$X$ measurable, implying that~$Y$ is independent of $Z$ conditionally on $X$.
\end{proof}

\begin{proof}[Proof of Theorem~\ref{thm:general-norms}]
    \ref{it:m1}: Since~$\gamma$ maps into the interval~${\lbrack0,1\rbrack}$, we only need to check that the quotient of the norms remains inside its domain. Non-negativity is obvious. As for the upper bound, it is easy to see that~${Q(Y)=Q(X,Y)}$ and therefore, by Lemma~\ref{lem:conditional} with~${Z=Y}$, the denominator in~\eqref{eq:norm-measure} is larger or equal than the numerator.
    \par \ref{it:m2}: The covariance matrix~${\Var(Q(X))}$ is the zero matrix if and only if the conditional distribution of~$Y$ given~$X$ does not actually depend on~$X$, meaning that~$X$ and~$Y$ are independent.

        \par \ref{it:m3}: If~$Y$ is $\sigma(X)$-measurable, then~${Q(X)=Q(Y)}$ and therefore~${\psi_{\|\cdot\|,\gamma}=1}$. On the other hand, if~$Y$ is not $\sigma(X)$-measurable, then~$Y$ is not independent of itself given~$X$ and therefore by Lemma~\ref{lem:conditional}
    \begin{equation*}
        \|\Var(Q(X))\|<\|\Var(Q(X,Y))\|=\|\Var(Q(Y))\|.
    \end{equation*}

    \par\ref{it:m4} and~\ref{it:m5}: These are direct consequences of Lemma~\ref{lem:conditional} and the monotonicity of~$\gamma$.

    \par \ref{it:m1-conditional}: By~\ref{it:m1} we have~${\psi_{\|\cdot\|,\gamma}((X,Z),Y)\leq 1}$. Together with~\ref{it:m4}, this implies the statement.

    \par \ref{it:m2-conditional}: This follows directly from~\ref{it:m5} and monotonicity of~$\gamma$.

    \par\ref{it:m3-conditional}: By~\ref{it:m3}, numerator and denominator are equal if and only if~$Y$ is~$\sigma(X,Z)$-measurable.
\end{proof}

\section*{Acknowledgments}
The authors thank Maximilian Ofner and Julius Baumhakel for helpful discussions.

\section*{Funding}
This research was funded in whole, or in part, by the Austrian Science Fund (FWF) [\href{https://doi.org/10.55776/P35520}{10.55776/P35520}]. For the purpose of open access, the authors have applied a CC~BY public copyright licence to any Author Accepted Manuscript version arising from this submission.

\bibliography{references2}

\begin{thebibliography}{36}
\providecommand{\natexlab}[1]{#1}

\bibitem[{Anderson(1962)}]{anderson_distribution_1962}
T.~W. Anderson, 1962.
\newblock On the {Distribution} of the {Two}-{Sample} {Cramer}-von {Mises}
  {Criterion}.
\newblock \emph{The Annals of Mathematical Statistics}, 33(3):1148--1159.
\newblock Publisher: Institute of Mathematical Statistics.

\bibitem[{Ansari and Fuchs(2025)}]{ansari_direct_2025}
J.~Ansari and S.~Fuchs, 2025.
\newblock A direct extension of {Azadkia} \& {Chatterjee}'s rank correlation to
  multi-response vectors.
\newblock Preprint.

\bibitem[{Azadkia and Chatterjee(2021)}]{azadkia_simple_2021}
M.~Azadkia and S.~Chatterjee, 2021.
\newblock A simple measure of conditional dependence.
\newblock \emph{The Annals of Statistics}, 49(6):3070--3102.

\bibitem[{Azadkia and Roudaki(2025)}]{azadkia_new_2025}
M.~Azadkia and P.~Roudaki, 2025.
\newblock A new measure of dependence: {Integrated} ${R}^2$.
\newblock Preprint.

\bibitem[{Balogoun et~al.(2022)Balogoun, Nkiet, and
  Ogouyandjou}]{balogoun_k-sample_2022}
A.~S.~K. Balogoun, G.~M. Nkiet, and C.~Ogouyandjou, 2022.
\newblock A k-{Sample} {Test} for {Functional} {Data} {Based} on {Generalized}
  {Maximum} {Mean} {Discrepancy}.
\newblock \emph{Lithuanian Mathematical Journal}, 62(3):289--303.

\bibitem[{Berrett et~al.(2020)Berrett, Wang, Barber, and
  Samworth}]{berrett_conditional_2020}
T.~Berrett, Y.~Wang, R.~Barber, and R.~Samworth, 2020.
\newblock The {Conditional} {Permutation} {Test} for {Independence} {While}
  {Controlling} for {Confounders}.
\newblock \emph{Journal of the Royal Statistical Society Series B: Statistical
  Methodology}, 82(1):175--197.

\bibitem[{Candès et~al.(2018)Candès, Fan, Janson, and
  Lv}]{candes_panning_2018}
E.~Candès, Y.~Fan, L.~Janson, and J.~Lv, 2018.
\newblock Panning for {Gold}: ‘{Model}-{X}’ {Knockoffs} for {High}
  {Dimensional} {Controlled} {Variable} {Selection}.
\newblock \emph{Journal of the Royal Statistical Society Series B: Statistical
  Methodology}, 80(3):551--577.

\bibitem[{Candès and Tao(2007)}]{candes_dantzig_2007}
E.~Candès and T.~Tao, 2007.
\newblock The {Dantzig} selector: {Statistical} estimation when $p$ is much
  larger than $n$.
\newblock \emph{The Annals of Statistics}, 35(6):2313--2351.

\bibitem[{Chatterjee(2021)}]{chatterjee_new_2021}
S.~Chatterjee, 2021.
\newblock A {New} {Coefficient} of {Correlation}.
\newblock \emph{Journal of the American Statistical Association},
  116(536):2009--2022.

\bibitem[{Chatterjee(2024)}]{chatterjee_survey_2024}
S.~Chatterjee, 2024.
\newblock A {Survey} of {Some} {Recent} {Developments} in {Measures} of
  {Association}.
\newblock In S.~Athreya, A.~G. Bhatt, and B.~V. Rao (editors),
  \emph{Probability and {Stochastic} {Processes}: {A} {Volume} in {Honour} of
  {Rajeeva} {L}. {Karandikar}}, pages 109--128. Springer Nature, Singapore.

\bibitem[{Chen and Guestrin(2016)}]{chen_xgboost_2016}
T.~Chen and C.~Guestrin, 2016.
\newblock {XGBoost}: {A} {Scalable} {Tree} {Boosting} {System}.
\newblock In \emph{Proceedings of the 22nd {ACM} {SIGKDD} {International}
  {Conference} on {Knowledge} {Discovery} and {Data} {Mining}}, pages 785--794.

\bibitem[{Combes(2023)}]{combes_extension_2023}
R.~Combes, 2023.
\newblock An extension of {McDiarmid}'s inequality.
\newblock Preprint.

\bibitem[{Cramér(1999)}]{cramer_mathematical_1999}
H.~Cramér, 1999.
\newblock \emph{Mathematical {Methods} of {Statistics}}.
\newblock Princeton University Press.

\bibitem[{Dang et~al.(2021)Dang, Nguyen, Chen, and Zhang}]{dang_new_2021}
X.~Dang, D.~Nguyen, Y.~Chen, and J.~Zhang, 2021.
\newblock A new {Gini} correlation between quantitative and qualitative
  variables.
\newblock \emph{Scandinavian Journal of Statistics}, 48(4):1314--1343.

\bibitem[{Deb et~al.(2020)Deb, Ghosal, and Sen}]{deb_measuring_2020}
N.~Deb, P.~Ghosal, and B.~Sen, 2020.
\newblock Measuring {Association} on {Topological} {Spaces} {Using} {Kernels}
  and {Geometric} {Graphs}.
\newblock Preprint.

\bibitem[{Dette and Kroll(2025)}]{dette_simple_2025}
H.~Dette and M.~Kroll, 2025.
\newblock A simple bootstrap for {Chatterjee}’s rank correlation.
\newblock \emph{Biometrika}, 112(1):asae045.

\bibitem[{Fuchs(2024)}]{fuchs_quantifying_2024}
S.~Fuchs, 2024.
\newblock Quantifying directed dependence via dimension reduction.
\newblock \emph{Journal of Multivariate Analysis}, 201:105266.

\bibitem[{Hall and Tajvidi(2002)}]{hall_permutation_2002}
P.~Hall and N.~Tajvidi, 2002.
\newblock Permutation tests for equality of distributions in high‐dimensional
  settings.
\newblock \emph{Biometrika}, 89(2):359--374.

\bibitem[{Heller et~al.(2012)Heller, Heller, and
  Gorfine}]{heller_consistent_2012}
R.~Heller, Y.~Heller, and M.~Gorfine, 2012.
\newblock A consistent multivariate test of association based on ranks of
  distances.
\newblock \emph{Biometrika}, 100(2):503--510.

\bibitem[{Henze(1987)}]{henze_fraction_1987}
N.~Henze, 1987.
\newblock On the {Fraction} of {Random} {Points} with {Specified}
  {Nearest}-{Neighbour} {Interrelations} and {Degree} of {Attraction}.
\newblock \emph{Advances in Applied Probability}, 19(4):873--895.

\bibitem[{Hörmann and Strenger(2025)}]{hormann_azadkiachatterjees_2024}
S.~Hörmann and D.~Strenger, 2025.
\newblock Azadkia–{Chatterjee}'s dependence coefficient for infinite
  dimensional data.
\newblock \emph{Bernoulli}.
\newblock Forthcoming.

\bibitem[{Kelly et~al.(2025)Kelly, Longhjohn, and Nottingham}]{kelly_uci_2025}
M.~Kelly, R.~Longhjohn, and K.~Nottingham, 2025.
\newblock The {UCI} {Machine} {Learning} {Repository}.
\newblock \href{https://archive.ics.uci.edu/}{https://archive.ics.uci.edu/}.

\bibitem[{Kolmogorov(1933)}]{kolmogorov_sulla_1933}
A.~Kolmogorov, 1933.
\newblock Sulla determinazione empirica di una legge di distribuzione.
\newblock \emph{Giornale dell'Instituto Italiano degli Attuari}, 4:83--91.

\bibitem[{Kroll(2025)}]{kroll_asymptotic_2025}
M.~Kroll, 2025.
\newblock Asymptotic {Normality} of {Chatterjee}'s {Rank} {Correlation}.
\newblock Preprint.

\bibitem[{Kruskal and Wallis(1952)}]{kruskal_use_1952}
W.~H. Kruskal and W.~A. Wallis, 1952.
\newblock Use of {Ranks} in {One}-{Criterion} {Variance} {Analysis}.
\newblock \emph{Journal of the American Statistical Association},
  47(260):583--621.
\newblock Publisher: ASA Website \_eprint:
  https://www.tandfonline.com/doi/pdf/10.1080/01621459.1952.10483441.

\bibitem[{LeCun et~al.(1998)LeCun, Bottou, Bengio, and
  Haffner}]{lecun_gradient-based_1998}
Y.~LeCun, L.~Bottou, Y.~Bengio, and P.~Haffner, 1998.
\newblock Gradient-based learning applied to document recognition.
\newblock \emph{Proceedings of the IEEE}, 86(11):2278--2324.

\bibitem[{Liu and Shang(2025)}]{liu_measuring_2025}
Y.~Liu and P.~Shang, 2025.
\newblock Measuring {Feature}-{Label} {Dependence} {Using} {Projection}
  {Correlation} {Statistic}.
\newblock Preprint.

\bibitem[{Pan et~al.(2020)Pan, Wang, Zhang, Zhu, and Zhu}]{pan_ball_2020}
W.~Pan, X.~Wang, H.~Zhang, H.~Zhu, and J.~Zhu, 2020.
\newblock Ball {Covariance}: {A} {Generic} {Measure} of {Dependence} in
  {Banach} {Space}.
\newblock \emph{Journal of the American Statistical Association},
  115(529):307--317.

\bibitem[{Rinott(1994)}]{rinott_normal_1994}
Y.~Rinott, 1994.
\newblock On normal approximation rates for certain sums of dependent random
  variables.
\newblock \emph{Journal of Computational and Applied Mathematics},
  55(2):135--143.

\bibitem[{Schilling(1986{\natexlab{a}})}]{schilling_multivariate_1986}
M.~Schilling, 1986{\natexlab{a}}.
\newblock Multivariate {Two}-{Sample} {Tests} {Based} on {Nearest} {Neighbors}.
\newblock \emph{Journal of the American Statistical Association},
  81(395):799--806.

\bibitem[{Schilling(1986{\natexlab{b}})}]{schilling_mutual_1986}
M.~Schilling, 1986{\natexlab{b}}.
\newblock Mutual and shared neighbor probabilities: finite- and
  infinite-dimensional results.
\newblock \emph{Advances in Applied Probability}, 18(2):388--405.

\bibitem[{{SORA}(2019)}]{sora_national_2019}
{SORA}, 2019.
\newblock National election results {Austria} 1919 - 2017 ({OA} edition).
\newblock
  \href{https://doi.org/10.11587/EQUDAL}{https://doi.org/10.11587/EQUDAL}.

\bibitem[{{Statistik Austria}(2025)}]{statistik_austria_statatlas_2025}
{Statistik Austria}, 2025.
\newblock {STATatlas} - {Bevölkerungsstand}.
\newblock
  \href{https://www.statistik.at/atlas/}{https://www.statistik.at/atlas/}.

\bibitem[{Strenger-Galvis(2025)}]{strenger-galvis_fdep_2025}
D.~Strenger-Galvis, 2025.
\newblock {FDEP}: {Investing} {Dependence} in {Functional} {Data}.
\newblock
  \href{https://github.com/danielstrenger/FDEP}{https://github.com/danielstrenger/FDEP}.

\bibitem[{Székely and Rizzo(2013)}]{szekely_distance_2013}
G.~Székely and M.~Rizzo, 2013.
\newblock The distance correlation $t$-test of independence in high dimension.
\newblock \emph{Journal of Multivariate Analysis}, 117:193--213.

\bibitem[{Székely et~al.(2007)Székely, Rizzo, and
  Bakirov}]{szekely_measuring_2007}
G.~Székely, M.~Rizzo, and N.~Bakirov, 2007.
\newblock Measuring and testing dependence by correlation of distances.
\newblock \emph{The Annals of Statistics}, 35(6):2769--2794.

\end{thebibliography}

\appendix
\renewcommand{\theLemma}{S\arabic{Lemma}} 
\setcounter{Lemma}{0} 
\renewcommand{\theequation}{S\arabic{equation}}
\setcounter{equation}{0} 
\section{Supplement: Proofs of auxiliary results}\label{sec:appendix}
In the sequel, we use the following notation: for a realization~${((x_i, y_i)\colon 1\leq i\leq n)}$ of the random variables~${((X_i, Y_i)\colon 1\leq i\leq n)}$ the realizations of~${L_n}$ and~${N(i)}$ are denoted by~${l_n}$ and~${n(i)}$. To indicate dependence on a particular~${\mathbf{x}=(x_1,\ldots, x_n)}$ or a nearest neighbour graph~$g$ we use~${n_\mathbf{x}(i)}$ and~${n_g(i)}$, respectively.
\subsection{Proof of Lemma~\ref{lem:convergence}}
For the proof of Lemma~\ref{lem:convergence}, we need the following three results:

\begin{Lemma}\label{lem:neighbour-convergence}
 Let~$X_i$ be  i.i.d.\ random elements in a separable metric space~$\mathcal{H}$. Then as~$n\to\infty$ we have (i)~$X_{N(1)}\convas X_1$  and (ii)~$P(N(1)=N(2))=o(1)$.
\end{Lemma}
For the proof, we refer to Lemma~1 and Lemma~3 in \cite{hormann_azadkiachatterjees_2024}.

\begin{Lemma}\label{lem:convergence-distribution-1}
Suppose Assumption~\ref{ass:convergence} holds. Then for~$n\to\infty$, we have
\begin{equation*}
    (Y_1,Y_{N(1)}, Y_2, Y_{N(2)})\convd (Y_1, Y_1', Y_2, Y_2').
\end{equation*}
\end{Lemma}
\begin{proof}
Since we are dealing with discrete data Assumption~\ref{ass:convergence}~(a)
is equivalent to
\begin{equation}
1\{f(X_{N(1)},u)=k\}\stackrel{P}{\to} 1\{f(X_1,u)=k\}\label{eq:altcond}
\end{equation} 
for almost all~$u$ and for all~$1\leq k\leq K$.
    By Part~(ii) of Lemma~\ref{lem:neighbour-convergence} we get
\begin{align*}
    &P(Y_1=k_1, Y_{N(1)} =l_1, Y_2=k_2, Y_{N(2)}=l_2)\\
    &\quad = P(Y_1=k_1, f(X_{N(1)}, U_1') = l_1, 
    Y_2=k_2, f(X_{N(2)}, U_2') = l_2)+o(1),
\end{align*}
where~$U_i'\stackrel{iid}{\sim} \mathcal U(0,1)$ and independent of all other variables. Hence,
\begin{align*}
    &P(Y_1=k_1, f(X_{N(1)}, U_1') = l_1, 
    Y_2=k_2, f(X_{N(2)}, U_2') = l_2)\\
    & = \int_{\lbrack 0,1\rbrack^2} \ev\left( 1\{Y_1=k_1, Y_2=k_2\} 1\{f(X_{N(1)}, u) = l_1\} 1\{f(X_{N(2)}, v) = l_2)\}\right)dudv.
\end{align*}
Due to \eqref{eq:altcond}, the integrand converges in probability to
$$
1\{Y_1=k_1, Y_2=k_2\} 1\{f(X_1, u) = l_1\} 1\{f(X_2, v) = l_2\}
$$
and since it is bounded, the limit can pushed inside integral and expected value, yielding the claimed weak convergence.

Alternatively, assume that Assumption~\ref{ass:convergence}~(b) holds.
    Let~${\mathbf{X}_n\vcentcolon=(X_1,\ldots, X_n)}$ and define~${\mathbf{x}=(x_1,\ldots,x_n)\in\mathcal{H}^n}$. We have that
    \begin{align}
        &P(Y_1=k_1, Y_{N(1)}=l_1, Y_2=k_2, Y_{N(2)}=l_2)\nonumber\\
        &\quad = \int_{\mathcal{H}^n} P\big(f(x_1,U_1)=k_1, f(x_{n(1)},U_{n(1)})=l_1,\nonumber\\ 
        &\qquad\qquad f(x_2,U_2)=k_2, f(x_{n(2)},U_{n(2)})=l_2\big)~dP_{\mathbf{X}_n}(\mathbf{x}).\label{eq:integral-prob}
    \end{align}
    We split the integral over regions
\begin{equation}
R_1\vcentcolon=\{\mathbf{x}\colon |\{1,2,n(1),n(2)\}|=4\}\quad \text{and}\quad R_2\vcentcolon=R_1^c.
\end{equation}
By the second part of Lemma~\ref{lem:neighbour-convergence}, it follows that
    \begin{equation*}
        P(\mathbf{X}_n\in R_2)\leq P(N(1)=2)+P(N(2)=1)+P(N(1)=N(2))=o(1).
    \end{equation*}
Since the~$U_i$ are independent of each other and of the~$X_i$, the integrand in~\eqref{eq:integral-prob} factorizes on~$R_1$, i.e.\ we may write    
\begin{align*}
        &P(Y_1=k_1, Y_{N(1)}=l_1, Y_2=k_2, Y_{N(2)}=l_2)\\
        &\quad = \int_{R_1} P(f(x_1,U_1)=k_1)P( f(x_{n(1)},U_{1}')=l_1)\\ &\qquad\quad P(f(x_2,U_2)=k_2)P(f(x_{n(2)},U_{2}')=l_2)~dP_{\mathbf{X}_n}(\mathbf{x})+o(1)\\
        &\quad = \int_{\mathcal H^n} Q_{k_1}(x_1)Q_{l_1}(x_{n(1)}) Q_{k_2}(x_2)Q_{l_2}(x_{n(2)})~dP_{\mathbf{X}_n}(\mathbf{x}) + o(1)\\
        &\quad = \int_{\mathcal H^n} Q_{k_1}(x_1)Q_{l_1}(x_{1}) Q_{k_2}(x_2)Q_{l_2}(x_{2})~dP_{\mathbf{X}_n}(\mathbf{x})\\
        &\qquad + \int_{\mathcal H^n} Q_{k_1}(x_1)(Q_{l_1}(x_{n(1)})-Q_{l_1}(x_{1})) Q_{k_2}(x_2)Q_{l_2}(x_{2})~dP_{\mathbf{X}_n}(\mathbf{x})\\
        &\qquad + \int_{\mathcal H^n} Q_{k_1}(x_1)Q_{l_1}(x_{n(1)}) Q_{k_2}(x_2)(Q_{l_2}(x_{n(2)})-Q_{l_2}(x_{2}))~dP_{\mathbf{X}_n}(\mathbf{x}) + o(1)
        \\
        &\quad =:T_{1,n}+T_{2,n}+T_{3,n}+o(1).
    \end{align*}
    Now we note that
    $T_{1,n}=P(Y_1=k_1,Y_1'=l_1,Y_2=k_2,Y_2'=l_2)$. Moreover, we have
\begin{align*}
    |T_{2,n}|&\leq \int_{\mathcal H^n}|Q_{l_1}(x_{n(1)})-Q_{l_1}(x_{1})| ~dP_{\mathbf{X}_n}(\mathbf{x})=\ev|Q_{\ell_1}(X_{N(1)})-Q_{\ell_1}(X_{1})|.
\end{align*}
By Assumption~\ref{ass:convergence}~(b) we have that $Q_{\ell_1}(x)$ is $P_X$-almost surely continuous, and hence by Lemma~\ref{lem:neighbour-convergence}~(i) it follows that $Q_{\ell_1}(X_{N(1)})\to Q_{\ell_1}(X_{1})$. By dominated convergence we conclude that $T_{2,n}\to 0$ and by the same arguments we also get~$T_{3,n}\to 0$.

\end{proof}

\begin{Lemma}[\cite{combes_extension_2023}]\label{lem:combes}
    Let~$\mathcal{X}_1,\ldots, \mathcal{X}_n$ be a collection of arbitrary metric spaces and set~${\mathcal{X}=\mathcal{X}_1\times\ldots\times\mathcal{X}_n}$. Let~$\mathcal{Y}\subset \mathcal{X}$ and let~$X$ be an~$\mathcal{X}$-valued random variable. Denote~${p=P(X\notin \mathcal Y)}$.
    If~$f:\mathcal{X}\to\R$ is such that for all pairs~$x,y\in\mathcal Y$ which are only distinct in the~$i$-th coordinate we have \begin{equation}\label{eq:bounded-differences}
        |f(x)-f(y)|\leq c_i,
    \end{equation}
    then it holds that \begin{equation}\label{eq:combes}
        P(|f(X)-\ev\left(f(X)|X\in\mathcal Y\right)|\geq \varepsilon)\leq 2p+2\exp\left(-\frac{2\left(\varepsilon-p\sum_{i=1}^n c_i\right))^2}{\sum_{i=1}^n c_i^2}\right).
    \end{equation}
\end{Lemma}
\begin{Remark} \cite{combes_extension_2023} states Lemma~\ref{lem:combes} for~$\mathcal{X}_i = \mathbb{R}^p$, but this is not relevant for the proof and it remains valid for any metric spaces.
\end{Remark}

\begin{proof}[Proof of Lemma~\ref{lem:convergence}]
    (i) We consider
    \begin{align*}
        \Var\left(\widehat p_{k,l}\right) &= \Var\left(\frac{1}{n}\sum_{i=1}^n 1\lbrace Y_i = k, Y_{N(i)}=l\rbrace\right)\\
         &=  \frac{1}{n^2}\sum_{i=1}^n\Var\left(1\lbrace Y_i = k, Y_{N(i)}=l\rbrace\right)\\
        &\qquad + \frac{1}{n^2}\sum_{i\neq j}^n \Cov\left(1\lbrace Y_i = k, Y_{N(i)}=l\rbrace, 1\lbrace Y_j = k, Y_{N(j)}=l\rbrace\right)\\
        &\leq \frac{1}{n}\Var\left(1\lbrace Y_1 = k, Y_{N(1)}=l\rbrace\right)\\
        &\qquad+|\Cov\left(1\lbrace Y_1 = k, Y_{N(1)}=l\rbrace, 1\lbrace Y_2 = k, Y_{N(2)}=l\rbrace\right)|.
    \end{align*}
    The first term is~$O(1/n)$. The second term is
    \begin{align*}
        &|P\left( Y_1 = k, Y_{N(1)}=l, Y_2 = k, Y_{N(2)}=l \right)\\
        &\qquad -P\left(  Y_1 = k, Y_{N(1)}=l \right)P\left(  Y_2 = k, Y_{N(2)}=l \right)|,
    \end{align*}
 and by Lemma~\ref{lem:convergence-distribution-1}  it converges to
  \begin{align*}
        &|P\left( Y_1 = k, Y_1'=l, Y_2 = k, Y_2'=l \right)\\
        &\qquad -P\left(  Y_1 = k, Y_1'=l \right)P\left(  Y_2 = k, Y_2'=l \right)|.
    \end{align*}
    By independence of~$(Y_1,Y_1')$ 
     and~$(Y_2,Y_2')$ the last expression is equal to zero. Likewise, it follows from Assumption~\ref{ass:convergence} that~$\ev  \widehat p_{k,l}\to p_{k,l}$, and thus (i) follows.
     
    \par (ii) We use Lemma~\ref{lem:combes} with~$\mathcal X=\mathcal{H}^n\times\{1,\ldots,K\}^n$,
    \begin{equation*}
        f(x_1,\ldots, x_n, y_1,\ldots, y_n)=\frac{1}{n}\sum_{j=1}^n 1\{y_j=k,y_{n(j)}=l\},
    \end{equation*}
    and
    \begin{equation*}
        \mathcal Y = \lbrace (x_1,\ldots, x_n,y_1,\ldots, y_n): l_n\leq k_n\rbrace.
    \end{equation*}
    We remark that~$\widehat p_{k,l}=f(X_1,\ldots, X_n,Y_1,\ldots, Y_n)$. 
    
    Consider some tuple~$(\mathbf{x},\mathbf{y})=(x_1,\ldots, x_n,y_1,\ldots, y_n)\in \mathcal Y$ and let~$(\mathbf{x}^{(i)},\mathbf{y})$ and~$(\mathbf{x},\mathbf{y}^{(i)})$ be perturbations of~$(\mathbf{x},\mathbf{y})$ where~$x_i$ is replaced by~$x_i'$ and~$y_i$ by~$y_i'$, while all other coordinates remain unchanged. Moreover, assume that~$(\mathbf{x}^{(i)},\mathbf{y})$ and~$(\mathbf{x},\mathbf{y}^{(j)})$ are also in~$\mathcal Y$. It follows that 
   ~$$
    |f(\mathbf{x},\mathbf{y})-f(\mathbf{x}^{(i)},\mathbf{y})|\leq \frac{1}{n}\sum_{j=1}^n 1\{n_\mathbf{x}(j)\neq n_{\mathbf{x}^{(i)}}(j)\}\leq \frac{2 k_n}{n},
   ~$$
    and
 ~$$
    |f(\mathbf{x},\mathbf{y})-f(\mathbf{x},\mathbf{y}^{(i)})|\leq \frac{1}{n}\sum_{j=1}^n (1\{j=i\}+1\{n(j)=i\})\leq \frac{1+ k_n}{n}.
   ~$$  
Therefore,~\eqref{eq:bounded-differences} holds with~${c_i=2k_n/n}$. Using~\eqref{eq:combes}, we obtain
    \begin{align*}
        & P\left(|\widehat p_{k,l}-\ev (\widehat p_{k,l}|(X_1,\ldots,X_n,Y_1,\ldots,Y_n)\in\mathcal Y_n)|\geq \varepsilon\right)\\
        &\quad\leq 2P(L_n\geq k_n) + 2\exp\left(-\frac{2(\varepsilon-2k_nP(L_n\geq k_n))^2}{4k_n^2/n}\right).
    \end{align*}
    By Assumption~\eqref{ref:assasconv}, this term is summable over~$n$ and hence by the Borel-Cantelli lemma it follows that
    \begin{equation*}
        \widehat p_{k,l}-\ev (\widehat p_{k,l}|(X_1,\ldots,X_n,Y_1,\ldots,Y_n)\in\mathcal Y_n)\convas 0.
    \end{equation*}
Likewise, using dominated convergence and again the Borel-Cantelli lemma we get
    \begin{equation*}
        \ev (\widehat p_{k,l}|(X_1,\ldots,X_n,Y_1,\ldots,Y_n)\in\mathcal Y_n)-\ev(\widehat p_{k,l})\convas 0.
    \end{equation*}
    As already shown above, it holds that~${\ev \widehat p_{k,l}\to p_{k,l}}$, which completes the proof.
\end{proof}

\begin{proof}[Proof of Theorem~\ref{thm:convergence} and Theorem~\ref{thm:convergenceas}]
Convergence of~$\widehat p_k$ and~$\widehat q_l$, respectively,  follows from Lemma~\ref{lem:convergence} as~$\widehat p_k=\sum_{l=1}^K\widehat p_{k,l}$ and~$\widehat q_k=\sum_{k=1}^K\widehat p_{k,l}$. The proof for both theorems is then immediate from the continuous mapping theorem.
\end{proof}


\subsection{Proof of Lemma~\ref{lem:clt}}

\begin{Lemma}\label{lem:determinant}
    It holds that
    \begin{align*}    \mathrm{det}\left(\Sigma_{n}\right)&=\left(1+W_{n}\right)^{K(K-1)/2}\left(1-W_{n}\right)^{(K-1)(K-2)/2}\prod_{i=1}^{K} p_i^{2(K-1)}.
    \end{align*}
\end{Lemma}
We defer the proof of Lemma~\ref{lem:determinant} to  Supplement~\ref{sec:calculation-sigma}.

\begin{Lemma}\label{thm:rinott}[\cite{rinott_normal_1994}, Thm. 2.2]
 Let~$V_1,\ldots,V_n$ be random variables having a dependence graph whose maximal degree is strictly less than~$k$ and satisfy
 \begin{enumerate}
 \item~$|V_i|\leq~B$ almost surely for~$i=1,\ldots,n$,
 \item~$\ev  V_i =0$ and
 \item~$\Var(\sum_{i=1}^n V_i)=1$.
 \end{enumerate}
 Then
 \begin{equation}\label{eq:rinott}
 \left|P\left(\sum_{i=1} V_i\leq~z\right)-\Phi(z)\right|\leq~\sqrt{\frac{1}{2\pi}}kB+16\sqrt{n}k^{3/2}B^2+10nk^2B^3.
 \end{equation}
\end{Lemma}

\begin{proof}[Proof of Lemma~\ref{lem:clt}]
Fix some~$g\in\mathcal{K}_n$ and set
\begin{align*}
V_{g,i}&\vcentcolon=\mathrm{vec}\Big((((1\lbrace Y_i =k\rbrace - p_k)(1\lbrace Y_{n_g(i)} =l\rbrace - p_l)))_{k,l=1}^{K-1}\Big).
    \end{align*}
   Observing that~$P_n = \frac{1}{\sqrt{n}}\sum V_{g,i}$ and using independence between $X$'s and $Y$'s we need to show  that
       \begin{equation*}
        \sup_{g\in\mathcal K_n}\sup_{z\in\R} \left|P\left(c'\frac{1}{\sqrt{n}}\sum_{i=1}^n\Sigma_{g,n}^{-1/2}V_{g,i}\leq z\right)-\Phi(z)\right|\to 0,\quad (n\to \infty).
    \end{equation*}
    Due to Lemma~\ref{lem:determinant} and the assumption~$W_n\leq 1-b/2$ there is a~$\delta>0$ such that~${\lambda_{\min}(\Sigma_{g,n})\geq \delta}$.
   Hence, $\Sigma_{g,n}^{-1/2}$ is well defined.
    We set~$V_i=c'\frac{1}{\sqrt{n}}\Sigma_{g,n}^{-1/2} V_{g,i}$ and apply Lemma~\ref{thm:rinott}. Conditions 2 and 3 are immediate from the definition of~$V_i$. Furthermore, notice that 
   ~$$
    |V_i|\leq \frac{\|c\|\| V_{g,i}\|\|\Sigma_{g,n}^{-1/2}\|}{\sqrt{n}}\leq \frac{1}{\sqrt{n\lambda_{\min}(\Sigma_{g,n})}}\leq \sqrt{\delta^{-1}}\ n^{-1/2}.
   ~$$

     Next we note that there is an edge between~$V_i$ and~$V_j$ in their dependency graph if and only if ~$N(i)=j$, or
        ~$N(j)=i$, or
       ~$N(i)=N(j)$. Therefore, the maximal degree of the dependence graph is of the same order as~$\mathrm{deg}(g)$. The condition~$\sup_{g\in\mathcal K_n}\mathrm{deg}(g)=o\left(n^{1/4}\right)$ thus implies that the right-hand side of~\eqref{eq:rinott} tends to zero as~$n\to\infty$.
\end{proof}

\section{Supplement: Calculation of~$\Sigma_n$}\label{sec:calculation-sigma}
In the sequel we calculate 
\begin{align*}
    &\ev\big((1\lbrace Y_i=k_1\rbrace -p_{k_1})(1\lbrace Y_{N(i)}=l_1\rbrace -p_{l_1})\times\\
        &\qquad\qquad
        (1\lbrace Y_j=k_2\rbrace -p_{k_2})(1\lbrace Y_{N(j)}=l_2\rbrace -p_{l_2})|\mathcal{G}_n\big).
\end{align*}
Summation over $1\leq i,j\leq n$ then provides the elements of $\Sigma_n$. 
The product above decomposes into 16 terms. Some of the terms  depend on the combination of~$i$ and~$j$.   Below we list all cases along with the number of occurrences. Setting
   ~$$
    W_n'=\frac{1}{n}\sum_{\substack{i,j=1\\ i\neq j}}^n 1\{N(i)=N(j)\}
   ~$$
we obtain:

\begin{enumerate}[left=0pt,itemsep=1ex,topsep=0pt]
    \item
    \begin{tabularx}{\linewidth}{@{}Xr@{}}
       ~$p_{k_1}p_{l_1}p_{k_2}p_{l_2}$ &~$[\#=n^2]$ \\
    \end{tabularx}

    \item
    \begin{tabularx}{\linewidth}{@{}Xr@{}}
       ~$-p_{l_1}p_{k_2}p_{l_2} P(Y_i=k_1|\mathcal{G}_n) = -p_{k_1}p_{l_1}p_{k_2}p_{l_2}$ &~$[\#=n^2]$ \\
    \end{tabularx}

    \item
    \begin{tabularx}{\linewidth}{@{}Xr@{}}
       ~$-p_{k_1}p_{l_1}p_{l_2} P(Y_j=k_2|\mathcal{G}_n) = -p_{k_1}p_{l_1}p_{k_2}p_{l_2}$ &~$[\#=n^2]$ \\
    \end{tabularx}

   \item
    \begin{tabularx}{\linewidth}{@{}Xr@{}}
       ~$p_{l_1}p_{l_2}P(Y_i=k_1, Y_j=k_2|\mathcal{G}_n)=$ & \\[1ex]
       ~$\quad\left\{
        \begin{array}{ll}
            1\lbrace k_1=k_2\rbrace p_{k_1}p_{l_1}p_{l_2} & \text{if } i = j \\
            p_{k_1}p_{l_1}p_{k_2}p_{l_2} & \text{otherwise}
        \end{array}
        \right.$ &
       ~$\begin{aligned}
            [\#=n] \\
            [\#=n^2 - n]
        \end{aligned}$
    \end{tabularx}

        \item
    \begin{tabularx}{\linewidth}{@{}Xr@{}}
       ~$-p_{k_1}p_{k_2}p_{l_2} P(Y_{N(i)}=l_1|\mathcal{G}_n) = -p_{k_1}p_{l_1}p_{k_2}p_{l_2}$ &~$[\#=n^2]$ \\
    \end{tabularx}

    \item
    \begin{tabularx}{\linewidth}{@{}Xr@{}}
       ~$p_{k_2}p_{l_2}P(Y_i=k_1, Y_{N(i)}=l_1|\mathcal{G}_n) = p_{k_1}p_{l_1}p_{k_2}p_{l_2}$ &~$[\#=n^2]$ \\
    \end{tabularx}

  \item
    \begin{tabularx}{\linewidth}{@{}Xr@{}}
       ~$p_{k_1}p_{l_2}P(Y_{N(i)}=l_1, Y_j=k_2|\mathcal{G}_n)=$ & \\[1ex]
        \quad~$\left\{
        \begin{array}{ll}
            1\lbrace l_1 = k_2 \rbrace p_{k_1}p_{l_1}p_{l_2} & \text{if } N(i) = j \\
            p_{k_1}p_{l_1}p_{k_2}p_{l_2} & \text{otherwise}
        \end{array}
        \right.$ &
       ~$\begin{aligned}
            [\#=n] \\
            [\#=n^2 - n]
        \end{aligned}$
    \end{tabularx}

        \item
    \begin{tabularx}{\linewidth}{@{}Xr@{}}
       ~$-p_{l_2}P(Y_i=k_1, Y_{N(i)}=l_1, Y_j=k_2|\mathcal{G}_n)=$ & \\[1ex]
       ~$\quad\left\{
        \begin{array}{ll}
            -1\lbrace k_1=k_2\rbrace p_{k_1}p_{l_1}p_{l_2} & \text{if } i = j \\
            -1\lbrace l_1=k_2\rbrace p_{k_1}p_{l_1}p_{l_2} & \text{if } N(i) = j \\
            -p_{k_1}p_{l_1}p_{k_2}p_{l_2} & \text{otherwise}
        \end{array}
        \right.$ &
       ~$\begin{aligned}
            [\#=n] \\
            [\#=n] \\
            [\#=n^2 - 2n]
        \end{aligned}$
    \end{tabularx}

    \item
    \begin{tabularx}{\linewidth}{@{}Xr@{}}
       ~$-p_{k_1}p_{l_1}p_{k_2}P(N(j)=l_2|\mathcal{G}_n) = -p_{k_1}p_{l_1}p_{k_2}p_{l_2}$ &~$[\#=n^2]$ \\
    \end{tabularx}

    \item
    \begin{tabularx}{\linewidth}{@{}Xr@{}}
       ~$p_{l_1}p_{k_2}P(Y_i=k_1, Y_{N(j)}=l_2|\mathcal{G}_n)=$ & \\[1ex]
        \quad~$\left\{
        \begin{array}{ll}
            1\lbrace k_1=l_2\rbrace p_{k_1}p_{l_1}p_{k_2} & \text{if } i = N(j) \\
            p_{k_1}p_{l_1}p_{k_2}p_{l_2} & \text{otherwise}
        \end{array}
        \right.$ &
       ~$\begin{aligned}
            [\#=n] \\
            [\#=n^2 - n]
        \end{aligned}$
    \end{tabularx}

    \item
    \begin{tabularx}{\linewidth}{@{}Xr@{}}
       ~$p_{k_1}p_{l_1}P(Y_j=k_2, Y_{N(j)}=l_2|\mathcal{G}_n) = p_{k_1}p_{l_1}p_{k_2}p_{l_2}$ &~$[\#=n^2]$ \\
    \end{tabularx}

    \item
    \begin{tabularx}{\linewidth}{@{}Xr@{}}
       ~$-p_{l_1}P(Y_i=k_1, Y_j = k_2, Y_{N(j)}= l_2|\mathcal{G}_n)=$ & \\[1ex]
       ~$\quad\left\{
        \begin{array}{ll}
            -1\lbrace k_1=k_2\rbrace p_{k_1}p_{l_1}p_{l_2} & \text{if } i = j \\
            -1\lbrace k_1=l_2\rbrace p_{k_1}p_{l_1}p_{k_2} & \text{if } i = N(j) \\
            -p_{k_1}p_{l_1}p_{k_2}p_{l_2} & \text{otherwise}
        \end{array}
        \right.$ &
       ~$\begin{aligned}
            [\#=n] \\
            [\#=n] \\
            [\#=n^2 - 2n]
        \end{aligned}$
    \end{tabularx}

    \item
    \begin{tabularx}{\linewidth}{@{}Xr@{}}
       ~$p_{k_1}p_{k_2}P(Y_{N(i)}=l_1, Y_{N(j)}=l_2|\mathcal{G}_n)=~$ & \\[1ex]
       ~$\quad\left\{
        \begin{array}{ll}
            1\lbrace l_1 = l_2 \rbrace p_{k_1}p_{l_1}p_{k_2} & \text{if } N(i) = N(j) \\
            p_{k_1}p_{l_1}p_{k_2}p_{l_2} & \text{otherwise}
        \end{array}
        \right.$ &
       ~$\begin{aligned}
            [\#=n(W_{n}' + 1)] \\
            [\#=n(n - 1 - W_{n}')]
        \end{aligned}$
    \end{tabularx}

    \item
    \begin{tabularx}{\linewidth}{@{}Xr@{}}
       ~$-p_{k_2}P(Y_i=k_1, Y_{N(i)}=l_1, Y_{N(j)}=l_2|\mathcal{G}_n)=~$ & \\[1ex]
       ~$\quad=\left\{
        \begin{array}{ll}
            -1\lbrace k_1=l_2\rbrace p_{k_1}p_{l_1}p_{k_2} & \text{if } i = N(j) \\
            -1\lbrace l_1=l_2\rbrace p_{k_1}p_{l_1}p_{k_2} & \text{if } N(i) = N(j) \\
            -p_{k_1}p_{l_1}p_{k_2}p_{l_2} & \text{otherwise}
        \end{array}
        \right.$ &
       ~$\begin{aligned}
            [\#=n] \\
            [\#=n(W_{n}' + 1)] \\
            [\#=n(n - 2 - W_{n}')]
        \end{aligned}$
    \end{tabularx}

    \item
    \begin{tabularx}{\linewidth}{@{}Xr@{}}
       ~$-p_{k_1}P(Y_{N(i)}=l_1, Y_j=k_2, Y_{N(j)}=l_2|\mathcal{G}_n)=~$ & \\[1ex]
       ~$\quad\left\{
        \begin{array}{ll}
            -1\lbrace l_1=k_2\rbrace p_{k_1}p_{l_1}p_{l_2} & \text{if } N(i) = j \\
            -1\lbrace l_1=l_2\rbrace p_{k_1}p_{l_1}p_{k_2} & \text{if } N(i) = N(j) \\
            -p_{k_1}p_{l_1}p_{k_2}p_{l_2} & \text{otherwise}
        \end{array}
        \right.$ &
       ~$\begin{aligned}
            [\#=n] \\
            [\#=n(W_{n}' + 1)] \\
            [\#=n(n - 2 - W_{n}')]
        \end{aligned}$
    \end{tabularx}

    \item
    \begin{tabularx}{\linewidth}{@{}Xr@{}}
       ~$P(Y_i=k_1, Y_{N(i)}=l_1, Y_j= k_2, Y_{N(j)} = l_2|\mathcal{G}_n) =$ & \\[1ex]
       ~$\quad\left\{
        \begin{array}{ll}
            1\lbrace k_1=k_2, l_1=l_2 \rbrace p_{k_1}p_{l_1} & \text{if } i = j \\
            1\lbrace k_1=l_2, l_1=k_2 \rbrace p_{k_1}p_{l_1} & \text{if } N(i)=j, N(j)=i \\
            1\lbrace l_1=l_2 \rbrace p_{k_1}p_{l_1}p_{k_2} & \text{if } N(i)=N(j),\ i \ne j \\
            1\lbrace k_1=l_2 \rbrace p_{k_1}p_{l_1}p_{k_2} & \text{if } i=N(j),\ N(i)\ne j \\
            1\lbrace l_1=k_2 \rbrace p_{k_1}p_{l_1}p_{l_2} & \text{if } N(i)=j,\ i \ne N(j) \\
            p_{k_1}p_{l_1}p_{k_2}p_{l_2} & \text{otherwise}
        \end{array}
        \right.$ &
       ~$\begin{aligned}
            [\# =n] \\[-0.8ex]
            [\#= nW_{n}] \\[-0.8ex]
            [\#= nW_{n}'] \\[-0.8ex]
            [\#=n(1 - W_{n})] \\[-0.8ex]
            [\#=n(1 - W_{n})] \\[-0.8ex]
            [\#=h_n]
        \end{aligned}$
    \end{tabularx}    
    \end{enumerate}
    where $h_n=n(n - 3 + W_{n} - W_{n}')$.

    \begin{Remark}
Note that the term~$W_{n}'$ doesn't appear in~$\Sigma_n$. The corresponding summands cancel.
\end{Remark}

\subsection{Proof of Lemma~\ref{lem:determinant}}\label{sec:proof-determinant}

To compute the determinant of~$\Sigma_n$ we use the following factorization.
\begin{Lemma}
Fix~$K$ and let~${E\in \mathbb{R}^{(K-1)\times (K-1)}}$ be a matrix with all entries equal to 1 and~${I}$ be the identity matrix in~${\mathbb{R}^{(K-1)\times (K-1)}}$.  Let~$D=\mathrm{diag}(p_1,\ldots, p_{K-1})$ and denote~$A^{2\otimes}=A\otimes A$. Finally, define the matrix~${\mathcal{W}_n\in\mathbb{R}^{(K-1)^2\times (K-1)^2}}$ with
$$\mathcal{W}_n((k_1,l_1),(k_2,l_2))=1\lbrace k_1=k_2, l_1=l_2\rbrace + W_{n} 1\lbrace k_1=l_2, l_1=k_2\rbrace,
$$
for~$1\leq k_1,l_1,k_2,l_2\leq K-1$.
Then the following factorization holds:
\begin{equation}\label{eq:expansionSigman}
\Sigma_n=D^{\otimes 2}\cdot \mathcal{W}_n\cdot (I-E D)^{2\otimes}.
\end{equation}
\end{Lemma}

\begin{proof}
First we use basic rules for Kronecker-products to expand the right-hand factor:
\begin{align*}
&(I-E D)^{2\otimes}\\
&=I\otimes I-(E\otimes I) (D\otimes I)-(I\otimes E)(I\otimes D)+ (E\otimes E)(D\otimes D).
\end{align*}
The involved matrices have the following form: 
\begin{flalign*}
(E\otimes I) (D\otimes I)&=
\begin{blockarray}{[ccc|ccc|ccc]}
p_1 & 0 & 0    & p_2 & 0    & 0   & p_{K-1} & 0   & 0 \\
0 & \ddots & 0 & 0 & \ddots & 0   & 0   & \ddots & 0 \\
0 & 0 & p_1    & 0 & 0      & p_2 & 0 & 0   & p_{K-1} \\
\cmidrule(lr){1-9}
p_1 & 0 & 0    & p_2 & 0    & 0   & p_{K-1} & 0   & 0 \\
0 & \ddots & 0 & 0 & \ddots & 0   & 0   & \ddots & 0 \\
0 & 0 & p_1    & 0 & 0      & p_2 & 0 & 0   & p_{K-1} \\
\cmidrule(lr){1-9}
p_1 & 0 & 0    & p_2 & 0    & 0   & p_{K-1} & 0   & 0 \\
0 & \ddots & 0 & 0 & \ddots & 0   & 0   & \ddots & 0 \\
0 & 0 & p_1    & 0 & 0      & p_2 & 0 & 0   & p_{K-1} \\
\end{blockarray}\\
(I\otimes E)(I\otimes D)&=
\begin{blockarray}{[ccc|ccc|ccc]}
p_1 & \cdots & p_{K-1} & 0 & \cdots & 0 & 0 & \cdots & 0 \\
\vdots & \ddots & \vdots & \vdots & \ddots & \vdots & \vdots & \ddots & \vdots \\
p_1 & \cdots & p_{K-1} & 0 & \cdots & 0 & 0 & \cdots & 0 \\
\cmidrule(lr){1-9}
0 & \cdots & 0 & p_1 & \cdots & p_{K-1} & 0 & \cdots & 0 \\
\vdots & \ddots & \vdots & \vdots & \ddots & \vdots & \vdots & \ddots & \vdots\\
0 & \cdots & 0 & p_1 & \cdots & p_{K-1} & 0 & \cdots & 0\\
\cmidrule(lr){1-9}
0 & \cdots & 0 & 0 & \cdots & 0 & p_1 & \cdots & p_{K-1} \\
\vdots & \ddots & \vdots & \vdots & \ddots & \vdots & \vdots & \ddots & \vdots \\
0 & \cdots & 0 & 0 & \cdots & 0 & p_1 & \cdots & p_{K-1} \\
\end{blockarray}
\end{flalign*}

\begin{flalign*}
(E\otimes E)(D\otimes D) &=
\begin{blockarray}{[ccc|ccc|ccc]}
p_1p_1 & \cdots & p_1p_{K-1} & p_2p_1 & \cdots & p_2p_{K-1} & p_{K-1}p_1 & \cdots & p_{K-1} p_{K-1} \\
\vdots & \ddots & \vdots & \vdots & \ddots & \vdots & \vdots & \ddots & \vdots \\
p_1p_1 & \cdots & p_1p_{K-1} & p_2p_1 & \cdots & p_2p_{K-1} & p_{K-1}p_1 & \cdots & p_{K-1} p_{K-1} \\
\cmidrule(lr){1-9}
p_1p_1 & \cdots & p_1p_{K-1} & p_2p_1 & \cdots & p_2p_{K-1} & p_{K-1}p_1 & \cdots & p_{K-1} p_{K-1} \\
\vdots & \ddots & \vdots & \vdots & \ddots & \vdots & \vdots & \ddots & \vdots \\
p_1p_1 & \cdots & p_1p_{K-1} & p_2p_1 & \cdots & p_2p_{K-1} & p_{K-1}p_1 & \cdots & p_{K-1} p_{K-1} \\
\cmidrule(lr){1-9}
p_1p_1 & \cdots & p_1p_{K-1} & p_2p_1 & \cdots & p_2p_{K-1} & p_{K-1}p_1 & \cdots & p_{K-1} p_{K-1} \\
\vdots & \ddots & \vdots & \vdots & \ddots & \vdots & \vdots & \ddots & \vdots \\
p_1p_1 & \cdots & p_1p_{K-1} & p_2p_1 & \cdots & p_2p_{K-1} & p_{K-1}p_1 & \cdots & p_{K-1} p_{K-1} \\
\end{blockarray}\\
\mathcal{W}_n &=
\begin{blockarray}{[ccc|ccc|cccc]}
1+W_n & 0 & 0 & 0 & 0 & 0 & 0 & 0 & 0 \\
0 & 1 & 0 & W_n & 0 & 0 & 0 & 0 & 0 \\
0 & 0 & 1 & 0 & 0 & 0 & W_n & 0 & 0 \\
\cmidrule(lr){1-9}
0 & W_n & 0 & 1 & 0 & 0 & 0 & 0 & 0 \\
0 & 0 & 0 & 0 & 1+W_n & 0 & 0 & 0 & 0 \\
0 & 0 & 0 & 0 & 0 & 1 & 0 & W_n & 0 \\
\cmidrule(lr){1-9}
0 & 0 & W_n & 0 & 0 & 0 & 1 & 0 & 0 \\
0 & 0 & 0 & 0 & 0 & W_n & 0 & 1 & 0 \\
0 & 0 & 0 & 0 & 0 & 0 & 0 & 0 & 1+W_n \\
\end{blockarray}
\end{flalign*}

Combining these matrices, \eqref{eq:expansionSigman} can be readily verified. 
\end{proof}

By basic properties of determinants we deduce from this lemma that 
$$
\mathrm{det}(\Sigma_n)=\mathrm{det}(D)^{2(K-1)}\cdot\mathrm{det}(\mathcal{W}_n)\cdot\mathrm{det}(I-E D)^{2(K-1)}.
$$
Clearly,~$\mathrm{det}(D)=\prod_{i=1}^{K-1}p_i$. The matrices~$\mathcal{W}_n$ and~$I-ED$ can be easily brought in upper diagonal forms:

\[
\mathrm{UD}(\mathcal{W}_n) =
\begin{blockarray}{[ccc|ccc|ccc]}
1+W_n & 0 & 0 & 0 & 0 & 0 & 0 & 0 & 0 \\
0 & 1 & 0 & W_n & 0 & 0 & 0 & 0 & 0 \\
0 & 0 & 1 & 0 & 0 & 0 & W_n & 0 & 0 \\
\cmidrule(lr){1-9}
0 & 0 & 0 & 1-W_n^2 & 0 & 0 & 0 & 0 & 0 \\
0 & 0 & 0 & 0 & 1+W_n & 0 & 0 & 0 & 0 \\
0 & 0 & 0 & 0 & 0 & 1 & 0 & W_n & 0 \\
\cmidrule(lr){1-9}
0 & 0 & 0 & 0 & 0 & 0 & 1-W_n^2 & 0 & 0 \\
0 & 0 & 0 & 0 & 0 & 0 & 0 & 1-W_n^2 & 0 \\
0 & 0 & 0 & 0 & 0 & 0 & 0 & 0 & 1+W_n \\
\end{blockarray}
\]
and
\[
\mathrm{UD}(I-ED) =
\begin{blockarray}{[ccc]}
1-p_1-\cdots -p_{K-1} & 0 & 0 \\
0 & 1 & 0  \\
0 & 0 & 1  \\
\end{blockarray}=
\begin{blockarray}{[ccc]}
p_K & 0 & 0 \\
0 & 1 & 0  \\
0 & 0 & 1  \\
\end{blockarray}.
\]
The proof of Lemma~\ref{lem:determinant} is now immediate.
\end{document}